\documentclass[10pt]{article}
\usepackage{amsmath,amsfonts,amscd,amssymb,hyperref}
\DeclareMathAlphabet\EuR{U}{eur}{m}{n}
\SetMathAlphabet\EuR{bold}{U}{eur}{b}{n}

\begin{document}

%%%%%%%%%%%%%%%%%%% Zaehler %%%%%%%%%%%%%%%%%%%%%%%%

%%\newtoks\theorembodyfont
%%\theorembodyfont{\slshape}
\newcommand\lp{\textup{(}}
\newcommand\rp{\textup{)}}
\newtheorem{theorem}{Theorem}[section]
\newtheorem{lemma}[theorem]{Lemma}
\newtheorem{notation}[theorem]{Notation}
\newtheorem{proposition}[theorem]{Proposition}
\newtheorem{definition}[theorem]{Definition}
\newtheorem{example}[theorem]{Example}
\newtheorem{remark}[theorem]{Remark}
\newtheorem{corollary}[theorem]{Corollary}
\newtheorem{conjecture}[theorem]{Conjecture}
\newtheorem{problem}[theorem]{Problem}
{\catcode`@=11\global\let\c@equation=\c@theorem}
\renewcommand{\theequation}{\thetheorem}
% Hier werden Gleichungen und Theoreme zusammen gezaehlt.
%Soll ein anderer Zaehler statt theorem verwendet werden
%(entspr. dem \newtheorem-Befehl), muss 2-mal theorem
%durch diesen Zaehler ersetzt werden. (Die Zeilen entsprechen
% der Zaehlung von \newtheorem{equation}[theorem]).

\renewcommand{\theenumi}{\alph{enumi}}
\renewcommand{\labelenumi}{\lp\theenumi\rp}
\renewcommand{\labelenumii}{\lp\theenumi\rp}

\providecommand{\proof}{\noindent{\bf \underline{\underline{Proof}} : }}
\providecommand{\qed}{\qquad \rule{2mm}{2mm} \bigskip}

\makeatletter
\renewcommand{\@seccntformat}[1]{\csname the#1\endcsname.\hspace{1em}}

\newcommand{\tit}[2]{\begin{bf} \begin{center} \begin{Large}
\section{#1}
\label{sec: #2}
\end{Large}\end{center}\end{bf}
\nopagebreak}

%%%%%%%%%%%%%%%%%%%%%%%%%%%%%% Diagrams %%%%%%%%%%%%%%%%%%%%%%%%%%%%%%

\newcommand{\squarematrix}[4]{\left( \begin{array}{cc} #1 & #2 \\ #3 &
#4
\end{array} \right)}
\newcommand{\smallmx}[4]{\mbox{\begin{scriptsize}$\squarematrix{#1}{#2}
        {#3}{#4}$\end{scriptsize}}}

\newcommand{\comsquare}[8]{
\begin{center}
$\begin{CD}
#1 @>#2>> #3\\
@V{#4}VV @VV{#5}V\\
#6 @>>#7> #8
\end{CD}$
\end{center}}

%%%%%%%%%%%%%%%%%% Categories %%%%%%%%%%%%%%%%%%%%%

\let\sect=\S
\newcommand{\curs}{\EuR}
\newcommand{\Chain}{\curs{Chain}}
\newcommand{\GROUPOIDS}{\curs{GROUPOIDS}}
\newcommand{\PAIRS}{\curs{PAIRS}}
\newcommand{\FGINJ}{\curs{FGINJ}}
\newcommand{\Mod}{\curs{Mod}}
\newcommand{\Or}{\curs{Or}}
\newcommand{\Sets}{\curs{Sets}}
\newcommand{\SPACES}{\curs{SPACES}}
\newcommand{\SPECTRA}{\curs{SPECTRA}}
\newcommand{\Sub}{\curs{Sub}}
\newcommand{\Orcat}{\Or(G;\calfin)}

%%%%%%%%%%%%%%% mathbb %%%%%%%%%%%%%%%%%%

\newcommand{\zz}{{\mathbb Z}}
\newcommand{\nn}{{\mathbb N}}
\newcommand{\cc}{{\mathbb C}}
\newcommand{\hh}{{\mathbb H}}
\newcommand{\kk}{{\mathbb K}}
\newcommand{\kko}{\mathbb{KO}}
\newcommand{\qq}{{\mathbb Q}}
\newcommand{\pp}{{\mathbb P}}
\newcommand{\rr}{{\mathbb R}}

\newcommand{\calb}{{\cal B}}
\newcommand{\calbh}{{\cal BH}}
\newcommand{\calc}{{\cal C}}
\newcommand{\cald}{{\cal D}}
\newcommand{\calf}{{\mathcal F}}
\newcommand{\calfin}{{\mathcal Fin}}
\newcommand{\calg}{{\cal G}}
\newcommand{\calh}{{\cal H}}
\newcommand{\calk}{{\cal K}}
\newcommand{\caln}{{\cal N}}
\newcommand{\cals}{{\cal S}}
\newcommand{\calz}{{\cal Z}}

%%%%%%%%%%%%%%%%%%%%% Operatornames %%%%%%%%%%%%%%%%%%%%%%%%%

\newcommand{\asmb}{\operatorname{asmb}}
\newcommand{\aut}{\operatorname{Aut}}
\newcommand{\ch}{\operatorname{ch}}
\newcommand{\chern}{\operatorname{chern}}
\newcommand{\CHAIN}{\operatorname{CHAIN}}
\newcommand{\class}{\operatorname{class}}
\newcommand{\clos}{\operatorname{clos}}
\newcommand{\cok}{\operatorname{coker}}
\newcommand{\cone}{\operatorname{cone}}
\newcommand{\colim}{\operatorname{colim}}
\newcommand{\con}{\operatorname{con}}
\newcommand{\conhom}{\operatorname{conhom}}
\newcommand{\consub}{\operatorname{consub}}
\newcommand{\cyclic}{\operatorname{cyclic}}
\newcommand{\Deg}{\operatorname{Deg}}
\newcommand{\diam}{\operatorname{diam}}
\newcommand{\ev}{\operatorname{ev}}
\newcommand{\Eul}{\operatorname{Eul}}
\newcommand{\Fix}{\operatorname{Fix}}
\newcommand{\Endo}{\operatorname{End}}
\newcommand{\Gen}{\operatorname{Gen}}
\newcommand{\HS}{\operatorname{HS}}
\newcommand{\hur}{\operatorname{hur}}
\newcommand{\inc}{\operatorname{inc}}
\newcommand{\id}{\operatorname{id}}
\newcommand{\ind}{\operatorname{ind}}
\newcommand{\im}{\operatorname{im}}
\newcommand{\inte}{\operatorname{int}}
\newcommand{\Irr}{\operatorname{Irr}}
\newcommand{\Is}{\operatorname{Is}}
\newcommand{\Isotr}{\operatorname{Isotr}}
\newcommand{\loc}{\operatorname{loc}}
\newcommand{\map}{\operatorname{map}}
\newcommand{\MOD}{\operatorname{MOD}}
\newcommand{\MODULES}{\operatorname{MODULES}}
\newcommand{\mor}{\operatorname{mor}}
\newcommand{\Ob}{\operatorname{Ob}}
\newcommand{\odd}{\operatorname{odd}}
\newcommand{\pr}{\operatorname{pr}}
\newcommand{\Rep}{\operatorname{Rep}}
\newcommand{\res}{\operatorname{res}}
\newcommand{\Sign}{\operatorname{Sign}}
\newcommand{\sign}{\operatorname{sign}}
\newcommand{\st}{\operatorname{st}}
\newcommand{\topo}{\operatorname{top}}
\newcommand{\Tor}{\operatorname{Tor}}
\newcommand{\tr}{\operatorname{tr}}
\newcommand{\vol}{\operatorname{vol}}
\newcommand{\Zero}{\operatorname{Zero}}

\newcommand{\bfE}{\ensuremath{\mathbf{E}}}

%%%%%%%%%%%%%%%%%%%%%%%%%% others %%%%%%%%%%%%%%%%%%%%%%%%%%%%%%%%%%%%%%

\newcommand{\comment}[1]                      %comment of the author
{
{{\bf Comment: } {\ttfamily #1}}
}

\newcommand{\tabtit}[2]
{
\ref{#2}. & #1
}

%%%%%%%%%%%%%%%%%%%%%%%%%  Begin of the text  %%%%%%%%%%%%%%%

\title{The equivariant Lefschetz fixed point theorem for 
proper cocompact $G$-manifolds}
\author{Wolfgang L\"uck and 
Jonathan Rosenberg\thanks{Partially supported by NSF grants 
DMS-9625336 and DMS-0103647.}}
\maketitle

%%%%%%%%%%%%%%%%%%%%%%%%%%%% Abstract  %%%%%%%%%%%%%%%%%%%%%%%%%%%%%%%%%%%
\typeout{-----------------------  Abstract  ------------------------}
\begin{abstract}
Suppose one is given a discrete group $G$, a cocompact proper
$G$-manifold $M$, and a $G$-self-map $f\colon M \to M$.
Then we introduce the equivariant Lefschetz class of $f$, which is
globally defined in terms of cellular chain complexes, and the  local
equivariant Lefschetz class of $f$,  
which is locally defined in terms of fixed point data. We prove the
equivariant Lefschetz 
fixed point theorem, which says that these two classes agree. As a
special case, we prove an equivariant Poincar\'e-Hopf Theorem,
computing the universal equivariant Euler characteristic in terms
of the zeros of an equivariant vector field, and also
obtain an orbifold Lefschetz fixed point theorem.
Finally, we prove a realization theorem for universal equivariant
Euler characteristics.
\smallskip


\end{abstract}

% %%%%%%%%%%%%%%%%%%%%%%%%%%%%%Introduction %%%%%%%%%%%%%%%%%%%%%%%%%%%%%%%%
\typeout{-----------------------  Introduction ------------------------}

\setcounter{section}{-1}
\tit{Introduction}{Introduction}

Let us recall the classical Lefschetz fixed point theorem. 
Let $f \colon M \to M$ be a smooth self-map of a compact smooth
manifold $M$, such that  $\Fix(f) \cap \partial M = \emptyset$ and for
each $x \in \Fix(f)$, the determinant 
of the linear map $(\id - T_xf)\colon  T_xM \to T_xM$ is different from
zero. Denote by $T_xM^c$ the one-point compactification of $T_xM$,
which is homeomorphic to a sphere. Let $(\id -T_xf)^c\colon T_xM^c \to
T_xM^c$ be the 
homeomorphism  induced by the self-homeo\-mor\-phism 
$(\id -T_xf)\colon T_xM \to T_xM$. Denote by $\deg((\id -T_xf)^c)$ its
degree, which is 
$1$ or $-1$, depending on whether $\det(\id -T_xf)$ is positive or negative.
Let 
$$L^{\zz[\{1\}]}(f) ~ := ~ \sum_{p \ge 0} (-1)^p \cdot
\tr_{\qq}(H_p(f;\qq)) ~ = ~ 
\sum_{p \ge 0} (-1)^p \cdot \tr_{\zz}(C_p(f))$$
be the classical Lefschetz number of $f$, where $H_p(f;\qq)$ is the
map on the singular 
homology with rational coefficients and $C_p(f)$ is the chain map on
the cellular $\zz$-chain complex 
induced by $f$ for some smooth triangulation of $M$. The Lefschetz
fixed point theorem says that 
under the conditions above the fixed point set $\Fix(f) = \{x \in M
\mid f(x) = x\}$ is finite and 
\begin{eqnarray}
L^{\zz[\{1\}]}(f) & = & \sum_{x \in \Fix(f)} \deg((\id -T_xf)^c).
\label{classical Lefschetz fixed point theorem}
\end{eqnarray}
For more information about it we refer for instance to
\cite{Brown(1971)}.

The purpose of this paper is to generalize this to the following equivariant
setting. Let $G$ be a (not necessarily finite) discrete group $G$.
A smooth $G$-manifold $M$ is a smooth manifold with an action of $G$ by
diffeomorphisms. It is called \emph{cocompact} if the quotient space
$G\backslash M$ is 
compact. It is \emph{proper} if the map $G\times M \to M\times M,~ 
(g,m)\mapsto (g\cdot m, m)$ is proper; when the action is cocompact,
this happens if and only if all isotropy
groups are finite. One can equip $M$ with the structure of a proper
finite $G$-$CW$-complex by an equivariant smooth triangulation
\cite{Illman(2000)}. The main result of this paper is

\begin{theorem}[Equivariant Lefschetz fixed point theorem] 
\label{the: equivariant Lefschetz theorem}
Let $G$ be a\linebreak
discrete group. Let $M$ be a cocompact proper
$G$-manifold {\lp}possibly with boundary{\rp} and 
let $f\colon  M \to  M$ be a
smooth $G$-map.  
Suppose that  $\Fix(f) \cap \partial M = \emptyset$ and for
each $x \in \Fix(f)$ the determinant 
of the linear map $\id - T_xf\colon  T_xM \to T_xM$ is different from
zero. 

Then $G\backslash\Fix(f)$ is finite, the  \emph{equivariant Lefschetz
class} of $f$  {\lp}see 
Definition \ref{def: equivariant Lefschetz class Lambda^G(f)}{\rp}
\begin{eqnarray*}
\Lambda^G(f) & \in & U^G(M)
\end{eqnarray*}
is defined in terms of cellular chain complexes, and 
the  \emph{local equivariant Lefschetz class} of $f$
{\lp}see Definition \ref{def: Lambda^G_{loc}(f)}{\rp}
\begin{eqnarray*}
\Lambda_{\loc}^G(f)
& \in & U^G(M)
\end{eqnarray*} 
is defined. Also $\Lambda^G(f)$ and $\Lambda_{\loc}^G(f)$ depend
only on the differentials $T_xf$ for $x \in \Fix(f)$, and
$$\Lambda^G(f) ~ = ~ \Lambda^G_{\loc}(f).$$
\end{theorem}

If $G$ is trivial, Theorem \ref{the: equivariant Lefschetz theorem}
reduces to \eqref{classical Lefschetz fixed point theorem}.
We emphasize that we want to treat arbitrary discrete groups and take
the component structure of the various fixed point sets into account.

In Section \ref{sec: The orbifold Lefschetz number}  we will define
the orbifold 
Lefschetz number, which can also be viewed as an $L^2$-Lefschetz number, 
and prove the orbifold Lefschetz fixed point theorem
\ref{the: the orbifold Lefschetz fixed point theorem} 
in Section \ref{sec: The orbifold Lefschetz fixed point theorem}. It
is both a key 
ingredient in the proof of and a special case of the equivariant
Lefschetz fixed point theorem \ref{the: equivariant Lefschetz theorem}. 

In Section
\ref{sec: The equivariant Lefschetz classes} we introduce the
\emph{equivariant Lefschetz class} 
$\Lambda^G(f)$, which is globally defined in terms of cellular chain
complexes, 
and in Section \ref{sec: The local equivariant Lefschetz class} we 
introduce the \emph{local equivariant Lefschetz class}
$\Lambda^G_{\loc}(f)$, which is locally defined in terms of the
differentials at the fixed points. 
These two are identified by the equivariant Lefschetz fixed point theorem 
\ref{the: equivariant Lefschetz theorem}, whose proof is completed in 
Section \ref{sec: The proof of the equivariant Lefschetz fixed point theorem}.

A classical result (the Poincar\'e-Hopf Theorem) says
that the Euler characteristic of a compact smooth manifold 
can be computed by counting (with signs) the zeros of a vector field which
is transverse to the 
zero-section and points outward at the boundary. This is a corollary of
the classical 
Lefschetz fixed point theorem \eqref{classical Lefschetz fixed point
theorem} via the associated flow.  
In Section \ref{sec: Euler characteristic and index of a vector field in the
equivariant setting} we will extend this result to
the equivariant setting for proper cocompact $G$-manifolds 
by defining the universal equivariant Euler characteristic, 
defining the index of an
equivariant vector field which is transverse to the
zero-section and points outward at the boundary,  and proving their
equality in 
Theorem \ref{the: chi^G(X) = i^G(v)}. As an illustration we explicitly
compute the universal equivariant Euler characteristic and
the local equivariant index of an equivariant vector field  for the
standard action of the infinite dihedral group on $\rr$ in 
Example \ref{exa: infinite dihedral group}.

To prove Theorem \ref{the: chi^G(X) = i^G(v)}  was one motivation for
this paper,
since it is a key ingredient in  \cite{Lueck-Rosenberg(2002b)}.
There a complete answer is given to the question of what information is
carried by the element $\Eul^G(M)\in KO^G_0(M)$, the class
defined by
the equivariant Euler operator for a 
proper cocompact $G$-manifold $M$. Rosenberg \cite{Rosenberg(1999a)}
has already settled this question in the non-equivariant case
by perturbing the Euler operator by a vector field
and using the classical result that the Euler characteristic can be
computed by counting 
the zeros of a vector field. The equivariant version of this strategy
will be applied in \cite{Lueck-Rosenberg(2002b)}, which requires
having Theorem \ref{the: chi^G(X) = i^G(v)} available.

In Section  \ref{sec: Constructing equivariant manifolds with given 
component structure and universal equivariant
Euler characteristic} we discuss the problem whether 
there exists a proper smooth $G$-manifold $M$
with prescribed sets $\pi_0(M^H)$ for $H \subseteq G$ 
such that $\chi^G(M)$ realizes a given element
in $U^G(M)$. A necessary and sufficient condition for this is given in
Theorem \ref{the: U^G realization theorem}. Again this will have
applications in \cite{Lueck-Rosenberg(2002b)}.

The paper is organized as follows:
\\[1mm]
\begin{tabular}{ll}
\ref{sec: The orbifold Lefschetz number}. & 
The orbifold Lefschetz number
\\
\ref{sec: The orbifold Lefschetz fixed point theorem}. &
The orbifold Lefschetz fixed point theorem
\\
\ref{sec: The equivariant Lefschetz classes}
& The equivariant Lefschetz classes
\\
\ref{sec: The local equivariant Lefschetz class}. & 
The local equivariant Lefschetz class
\\
\ref{sec: The proof of the equivariant Lefschetz fixed point theorem}.
& The proof of the equivariant Lefschetz fixed point theorem
\\
\ref{sec: Euler characteristic and index of a vector field in the
equivariant setting}. 
& Euler characteristic and index of a vector field in the equivariant setting
\\
\ref{sec: Constructing equivariant manifolds with given 
component structure and universal equivariant
Euler characteristic}. & 
Constructing equivariant manifolds with given component structure and 
\\ & universal equivariant Euler characteristic
\\
 & References\\[1mm]
\end{tabular}

%%%%%%%%%%%%%%%%%%%%%%%%%%   Section 1 %%%%%%%%%%%%%%%%%%%%%%%%%% %%%%%%%%%%%%

\typeout{-----------------------  Section 1  ------------------------}

\tit{The orbifold Lefschetz number}{The orbifold Lefschetz number}

In order to define the various Lefschetz classes and prove the various
Lefschetz fixed point 
theorems  for cocompact proper $G$-manifolds, we  need some input
about traces.  

Let $R$ be a commutative associative ring with unit,
for instance $R = \zz$ or $R = \qq$.
Let $u \colon P \to P$ be an endomorphism of a finitely generated
projective $RG$-module. Choose a finitely generated projective $RG$-module $Q$  
and an isomorphism $v \colon  P \oplus Q \xrightarrow{\cong}
\bigoplus_{i \in I} RG$ for some finite index set $I$. We obtain
an $RG$-endomorphism 
$$v  \circ (u \oplus 0) \circ v^{-1}\colon    \bigoplus_{i \in I} RG \to
\bigoplus_{i \in I} RG.$$
Let $A = (a_{i,j})_{i,j \in I}$ be the matrix
associated to this map, i.e.,
$$v  \circ (u \oplus 0) \circ v^{-1}(\{w_i \mid i \in I\}) ~ = ~
\left\{\left.\sum_{i \in J} w_i \cdot a_{i,j}~ \right| ~ j \in I\right\}.$$
Define
\begin{eqnarray} 
\tr_{RG} \colon RG & \to & R,  \hspace{5mm} \sum_{g \in G} r_g \cdot g ~ \mapsto ~
r_1
\label{tr_{qq G} : qqG to qq}
\end{eqnarray}
where $r_1$ is the coefficient of the unit element $1 \in
G$. Define the \emph{$RG$-trace of $u$} by 
\begin{eqnarray}
\tr_{RG}(u) & := & \sum_{i \in I} \tr_{RG}(a_{ii})
\hspace{5mm} \in R. \label{tr_{RG}(a) in R}
\end{eqnarray}

We omit the easy and well-known proof that this definition is independent of the various
choices such as $Q$ and $v $ and that the following Lemma 
\ref{lem: basic properties of tr_{qq G}(a)} is true.

\begin{lemma} \label{lem: basic properties of tr_{qq G}(a)}
\begin{enumerate}

\item  \label{lem: basic properties of tr_{qq G}(a): trace property}
Let $u \colon P \to Q$ and $v \colon Q \to P$ be $R G$-maps
of finitely generated projective $RG$-modules. Then
$$\tr_{RG}(v \circ u) = \tr_{RG}(u \circ v);$$

\item  \label{lem: basic properties of tr_{qq G}(a): diagonal}
Let $P_1$ and $P_2$ be finitely generated projective $RG$-modules. Let
$$\squarematrix{u_{1,1}}{u_{1,2}}{u_{2,1}}{u_{2,2}}\colon
P_1 \oplus P_2 \to P_1 \oplus P_2$$
be a $RG$-self-map. Then
$$\tr_{RG} \squarematrix{u_{1,1}}{u_{1,2}}{u_{2,1}}{u_{2,2}}
~ = ~ \tr_{RG}(u_{1,1}) + \tr_{RG}(u_{2,2});$$

\item \label{lem: basic properties of tr_{qq G}(a): linearity}
Let $u_1,u_2\colon P \to P$ be $RG$-endomorphisms of a finitely
generated projective $RG$-module and $r_1, r_2 \in R$. Then
$$\tr_{RG}(r_1 \cdot u_1  + r_2 \cdot u_2) 
~ = ~
r_1 \cdot \tr_{RG}(u_1) + r_2 \cdot \tr_{RG}(u_2);$$

\item  \label{lem: basic properties of tr_{qq G}(a): induction}
Let $\alpha\colon G \to K$ be an inclusion of groups and $u \colon P
\to P$ be an endomorphism 
of a finitely generated projective $RG$-module. Then induction with $\alpha$ 
yields an endomorphism $\alpha_*u$ of a finitely generated projective
$RK$-module, and 
$$\tr_{RK}(\alpha_*u) ~ = ~ \tr_{RG}(u);$$

\item  \label{lem: basic properties of tr_{qq G}(a): restriction}
Let $\alpha\colon H \to G$ be an inclusion of groups with finite index 
$[G:H]$ and $u \colon  P \to P$ be an endomorphism
of a finitely generated projective $RG$-module. Then the restriction to
$RH$ with $\alpha$ yields an  endomorphism $\alpha^*u$ 
of a finitely generated projective $RH$-module, and 
$$\tr_{RH}(\alpha^*u) ~ = ~ [G:H] \cdot \tr_{RG}(u);$$

\item  \label{lem: basic properties of tr_{qq G}(a): dimension and trace}
Let $H \subseteq G$ be finite such that $|H|$ is invertible in $R$.
Let $u\colon R[G/H] \to R[G/H]$ be
a $RG$-map which sends $1H$ to $\sum_{gH\in G/H} r_{gH} \cdot gH$. Then
$R[G/H]$ is a finitely  generated projective $RG$-module and 
\begin{eqnarray*}
\tr_{RG}(u) & = & |H|^{-1} \cdot r_{1H};
\\
\tr_{RG}(\id_{R[G/H]}) & = & |H|^{-1}.
\end{eqnarray*}

\end{enumerate}
\end{lemma}

Let $G$ be a discrete group. A relative $G$-$CW$-complex
$(X,A)$ is \emph{finite} if and only if  $X$ is obtained
from $A$ by attaching finitely many equivariant cells, or, equivalently, 
$G\backslash (X/A)$ is compact. A relative $G$-$CW$-complex $(X,A)$ is
\emph{proper} if and only if the isotropy group $G_x$ of each point $x
\in X-A$ is finite (see for instance \cite[Theorem 1.23 on page 18]{Lueck(1989)}).
Let $(f,f_0) \colon (X,A) \to (X,A)$ be a cellular $G$-self-map of a
finite proper relative $G$-$CW$-complex $(X,A)$.
Let $R$ be a commutative ring such that
for any $x \in X-A$ the order of its isotropy group $G_x$ is
invertible in $R$. 
Then the cellular $RG$-chain complex $C_*(X,A)$ is
finite projective, i.e., each chain module is finitely generated
projective and $C_p(X,A)= 0$ for $p \ge d$ for some integer $d$.

\begin{definition} \label{def: orbifold Lefschetz number}
Define the \emph{orbifold Lefschetz number} of $(f,f_0)$ by 
\begin{eqnarray}
L^{RG}(f,f_0) & := & \sum_{p \ge 0} (-1)^p \cdot \tr_{RG}(C_p(f,f_0))
\hspace{5mm} \in R.
\label{definition of L^{RG}(f)}
\end{eqnarray}
\end{definition}

One easily proves using Lemma \ref{lem: basic properties of tr_{qq G}(a)}

\begin{lemma} \label{lem: properties of L^{RG}(f)}
 Let $(f,f_0)\colon (X,A) \to (X,A)$ be a cellular $G$-self-map of
a finite proper relative $G$-$CW$-complex such that $|G_x|$ is
invertible in $R$ for 
each $x \in X-A$.  Then:

\begin{enumerate}
\item \label{lem: properties of L^{RG}(f): G-homotopy invariance}
The equivariant Lefschetz number  $L^{RG}(f,f_0)$ depends only on the
$G$-homo\-topy class of $(f,f_0)$; 

\item \label{lem: properties of L^{RG}(f): conjugation invariance}
Let $(g,g_0)\colon  (X,A) \to (Y,B)$ and $(h,h_0)\colon (Y,B) \to (X,A)$ be
cellular $G$-maps 
of finite proper relative $G$-$CW$-complexes such that 
$|G_x|$ is invertible in $R$ for each $x \in X-A$ and $|G_y|$ is
invertible in $R$ for each $y \in Y-B$. Then
$$L^{RG}(g \circ h,g_0 \circ h_0) ~ =  ~  L^{RG}(h \circ g,h_0 \circ g_0);$$

\item \label{lem: properties of L^{RG}(f): induction}
Let $\alpha\colon  G \to K$ be an inclusion of groups.
Then induction with $\alpha$ yields a cellular $K$-self-map
$\alpha_*(f,f_0)$ of a finite proper relative $K$-$CW$-complex, and 
$$L^{RK}(\alpha_*(f,f_0)) ~ = ~ L^{RG}(f,f_0);$$

\item \label{lem: properties of L^{RG}(f): restriction}
Let $\alpha\colon  H \to G$ be an inclusion of groups with finite index
$[G:H]$. 
Then restriction  with $\alpha$ yields a cellular $H$-self-map
$\alpha^*(f,f_0)$ of a finite proper relative $H$-$CW$-complex, and 
$$L^{RH}(\alpha^*f) ~ = ~ [G:H] \cdot L^{RG}(f).$$

\end{enumerate}
\end{lemma}

\begin{remark} \label{rem: L^G(f)}
\em The rational number $L^{\qq G}(f,f_0)$ agrees with the
\emph{$L^2$-Lefschetz number} 
$L^{(2)}(f,f_0;\caln(G))$ introduced in \cite[Section 6.8]{Lueck(2002)}.
It can be read off from the map induced by $(f,f_0)$ on the $L^2$-homology of 
$(X,X_0)$ by the analog of the usual formula, namely by 
$$L^{\qq G}(f,f_0) ~ = ~ L^{(2)}(f,f_0;\caln(G)) ~ = ~
\sum_{p \ge 0} (-1)^p \cdot \tr_{\caln(G)}\left(H_p^{(2)}(f,f_0;\caln(G))\right),$$
where $\tr_{\caln(G)}$ is the standard trace of the group von Neumann algebra $\caln(G)$. 
A similar formula exists in terms of $H_n(X,X_0;\qq)$ only under the very restrictive assumption,
that each $\qq G$-module $H_p(X,X_0;\qq)$ is finitely generated projective. 
If $G$ acts freely, then $L^{\qq G}(f)$ agrees with the (ordinary)
Lefschetz number $L^{\zz[\{1\}]}(G\backslash(f,f_0))$ of the
cellular self-map $G\backslash (f,f_0)$ of the finite relative
$CW$-complex $G\backslash (X,A)$. If $G$ is finite, then
$(X,X_0)$ is a finite relative $CW$-complex and 
$$L^{\qq G}(f,f_0) ~ = ~ \frac{1}{|G|} \cdot L^{\zz[\{1\}]}(f,f_0).$$
\em
\end{remark}

The following description of $L^G(f)$ will be useful later.
Let $I_p(X,A)$ be the set of path components of $X_p -X_{p-1}$. This
is the same as the set 
of open cells  of $(X,A)$ regarded as
relative $CW$-complex (after forgetting the group action).
The group $G$ acts on $I_p(X,A)$. For an open $p$-cell $e$ let $G_e$ be
its  isotropy group, 
$\overline{e}$ be its closure and $\partial e = e - \overline{e}$. Then
$\overline{e}/\partial e$ is homeomorphic to $S^p$ and there is a
homeomorphism 
$$h\colon  \bigvee_{e' \in I_p(X,A)} \overline{e'}/\partial e' ~
\xrightarrow{\cong} X_p/X_{p-1}.$$ 
For an open cell $e \in I_p(X,A)$ define the \emph{incidence number}
\begin{eqnarray} 
\inc(f,e) & \in & \zz
\label{definition of incidence number inc(f,e)}
\end{eqnarray} 
to be the degree of the composition
\begin{multline*}
\overline{e}/\partial e \xrightarrow{i_e} \bigvee_{e' \in I_p(X,A)} ~
 \overline{e'}/\partial e'  
 \xrightarrow{h } X_p/X_{p-1} 
\\
\xrightarrow{f} X_p/X_{p-1} \xrightarrow{h^{-1}} \bigvee_{e' \in I_p(X,A)} ~
 \overline{e'}/\partial e' 
\xrightarrow{\pr_e} \overline{e}/\partial e,
\end{multline*}
where $i_e³$ is the obvious inclusion and $\pr_e$ is the obvious projection.
Obviously $\inc(f,e) = \inc(f,ge)$ for $g \in G$. One easily checks
using Lemma \ref{lem: basic properties of tr_{qq G}(a)}

\begin{lemma} \label{lem: decription of L^{QG}(f) in terms of
incidence numbers} 
Let $(X,A)$ be a finite proper relative $G$-$CW$-complex. Consider a cellular
$G$-map $(f,f_0)\colon (X,A) \to (X,A)$. Then
$$L^{\qq G}(f,f_0) ~ = ~ \sum_{p \ge 0} (-1)^p \cdot 
\sum_{Ge \in G\backslash I_p(X,A)} |G_e|^{-1} \cdot \inc(f,e).$$
\end{lemma}

%%%%%%%%%%%%%%%%%%%%%%%%%%   Section 2 %%%%%%%%%%%%%%%%%%%%%%%%%% %%%%%%%%%%%%

\typeout{-----------------------  Section 2  ------------------------}

\tit{The orbifold Lefschetz fixed point theorem}
{The orbifold Lefschetz fixed point theorem}

This section is devoted to the proof of:
\begin{theorem}[The orbifold Lefschetz fixed point theorem] 
\label{the: the orbifold Lefschetz fixed point theorem}
Let $M$ be a cocompact proper $G$-manifold {\lp}possibly with
boundary{\rp} and let $f\colon  M \to 
M$ be a smooth $G$-map. Suppose that 
$\Fix(M) \cap \partial M = \emptyset$  and
for any $x \in Fix(f)$ the determinant of the map
$(\id_{T_xM} - T_{x}f)$ is different from zero. Then
$G\backslash \Fix(f)$ is finite, and 
$$L^{\qq G}(f) ~ = ~ \sum_{G\backslash \Fix(f)} |G_x|^{-1} \cdot 
\deg\left(\left(\id_{T_xM} - T_{x}f)\right)^c\right).$$ 
\end{theorem}

Theorem \ref{the: the orbifold Lefschetz fixed point theorem} above 
will be a key ingredient in the proof of
the equivariant Lefschetz fixed point theorem 
\ref{the: equivariant Lefschetz theorem}. On the other hand
Theorem \ref{the: equivariant Lefschetz theorem} implies 
Theorem \ref{the: the orbifold Lefschetz fixed point theorem}.

Let us first consider as an illustration the easy case, where $G$ is finite. Then
$$L^{\qq G}(f) = |G|^{-1} \cdot L^{\qq[\{1\}]}(f) = |G|^{-1} \cdot L^{\zz[\{1\}]}(f)$$
by Lemma \ref{lem: properties of L^{RG}(f)}
\eqref{lem: properties of L^{RG}(f): restriction} and
$L^{\zz[\{1\}]}(f)$ is the (ordinary) Lefschetz number of the self-map
$f\colon M \to M$ of the compact manifold $M$. The non-equivariant
Lefschetz fixed point theorem says
$$L^{\zz\{1\}}(f) ~ = ~ \sum_{\Fix(f)} 
\deg\left(\left(\id_{T_xM} - T_{x}f\right)^c\right).$$
Thus Theorem \ref{the: the orbifold Lefschetz fixed point theorem} 
follows for finite $G$. The proof in the case of an infinite group
cannot be reduced to the non-equivariant case in such an easy way 
since $M$ is not compact anymore. Instead we extend the proof in the
non-equivariant case to the equivariant setting.
\\
{\bf Proof of Theorem \ref{the: the orbifold Lefschetz fixed point theorem}:}
%For any $x \in \Fix(f)$ and $v \in T_xM$ with $T_xf(v) = v$ we have $v
%= 0$ by assumption. 
%From the argument below applied to a given  element $x \in \Fix(f)$ we get 
%a neighborhood $U_x$ such that $U_x \cap \Fix(f) = \{x\}$.
%Since $M$ is proper and cocompact,  $G\backslash\Fix(f)$ is finite.

Fix a $G$-invariant Riemannian metric on $M$. Choose $\epsilon_1 > 0$
such that for all $x \in M$ the exponential map is defined on 
$D_{\epsilon_1}T_xM = \{v \in T_xM \mid \Vert v\Vert \le \epsilon_1\}$,
where $\Vert v\Vert$ for $v \in T_xM$ is the norm coming from the Riemannian
metric. Such $\epsilon_1 > 0$ exists because $G\backslash M$ is compact.
The image  $N_{x,\epsilon_1}$ of the exponential map 
on $D_{\epsilon_1}T_xM$ is a $G_x$-submanifold of $M$ and a compact 
neighborhood of $x$. The exponential map induces a $G_x$-diffeomorphism
$$\exp_{x,\epsilon_1}\colon D_{\epsilon_1}T_xM \xrightarrow{\cong}
N_{x,\epsilon_1}$$ 
with $\exp_{x,\epsilon_1}(0) = x$ whose differential at $0$ is the 
identity  under the canonical identification
$T_0D_{\epsilon_1}T_xM = T_xM$. 

Since $G\backslash M$ is compact, we can 
choose $\epsilon_2 > 0$ such that $f(N_{x,\epsilon_2}) \subseteq
N_{x,\epsilon_1}$ and  $T_xf(D_{\epsilon_2}T_xM) \subseteq D_{\epsilon_1}T_xM$ holds for all
$x \in \Fix(f)$.  
Notice that $\exp_{x,\epsilon_1}$ restricted to $D_{\epsilon_2}T_xM$ is
$\exp_{x,\epsilon_2}$. We want to change
$f$ up to $G$-homotopy without changing $\Fix(f)$ 
such that $\exp_{x,\epsilon_1}^{-1} \circ f \circ \exp_{x,\epsilon_2}$ 
and $T_xf$ agree on $D_{\epsilon_3}T_xM$ for some positive number
$\epsilon_3 > 0$ and all $x \in \Fix(f)$.  Consider $x \in \Fix(f)$.
Notice that $\exp_{x,\epsilon_1}^{-1} \circ f \circ \exp_{x,\epsilon_2}$ sends
$0$ to $0$ and has $T_xf$ as differential at $0$ under the canonical
identification 
$T_0D_{\epsilon_1}T_xM = T_xM$. By Taylor's theorem
we can find a constant $C_1 > 0$ such that with respect to the norm on
$T_xM$ induced 
by the Riemannian metric on $M$
\begin{eqnarray}
\hspace{-5mm} 
||\exp_{x,\epsilon_1}^{-1} \circ f \circ \exp_{x,\epsilon_2}(v) - T_xf(v) || 
& \le &
C_1 \cdot ||v||^2 \hspace{2mm} \mbox{ for } v \in D_{\epsilon_2}T_xM.
\label{estimate efe^{-1}(v) - T_xf(v)}
\end{eqnarray} 
Since $\det(\id - T_xf) \not= 0$, we can find a constant $C_2 > 0$ such that
\begin{eqnarray}
||T_xf(v) - v || & \ge &
C_2 \cdot ||v|| \hspace{5mm} \mbox{ for }v \in T_xM.
\label{||T_xf(v) - v || ge C_2 cdot ||v||}
\end{eqnarray} 
Choose a smooth function $\phi\colon [0,\epsilon_2] \to [0,1]$ with
the properties that 
$\phi(t) = 1$ for $t \le \min\{C_2/3C_1,\epsilon_2/3\}$ and $\phi(t) = 0$ for 
$t \ge \min\{C_2/2C_1,\epsilon_2/2\}$. Define
$$h \colon D_{\epsilon_2}T_xM \times [0,1]  ~ \to ~ D_{\epsilon_1 }T_xM$$
by 
$$h(y,t) ~ := ~  \bigl(1 - t\phi(||v||)\bigr) \cdot  
\exp_{x,\epsilon_1}^{-1} \circ f \circ \exp_{x,\epsilon_2}(v)
 + t\phi(||v||) \cdot
T_xf(v).$$
Obviously $h$ is a $G_x$-homotopy from $h_0 = \exp_{x,\epsilon_1}^{-1}
\circ f \circ 
\exp_{x,\epsilon_2}$ to a $G_x$-map $h_1$. The homotopy $h$ is
stationary outside 
$D_{\min\{C_2/2C_1,\epsilon_2/2\}}T_xM$ and $h_1$ agrees with 
$T_xf$ on $D_{\epsilon_3}T_xM$ if we put $\epsilon_3 = \min\{C_2/3C_1,\epsilon_2/3\}$. 
Each map $h_t$ has on  $D_{\min\{C_2/2C_1,\epsilon_2/2\}}T_xM$ only one
fixed point, namely $0$. This follows from the following estimate 
based on \eqref{estimate efe^{-1}(v) - T_xf(v)} and 
\eqref{||T_xf(v) - v || ge C_2 cdot ||v||} for 
$v \in D_{\min\{C_2/2C_1,\epsilon_2/2\}}T_xM$:
\begin{eqnarray*}
\lefteqn{||h_t(v) - v||} & & 
\\ & \ge & ||T_xf(v) -v|| - 
\\ & & \hspace{10mm}||(1 - t\phi(||v||)) \cdot  
\exp_{x,\epsilon_1}^{-1} \circ f \circ \exp_{x,\epsilon_2}(v) -
(1 - t\phi(||v||)) \cdot  T_xf(v)||
\\
&  = &
||T_xf(v) -v|| - (1 - t\phi(||v||)) \cdot  || 
\exp_{x,\epsilon_1}^{-1} \circ f \circ \exp_{x,\epsilon_2}(v) - T_xf(v)||
\\
& \ge & 
C_2 \cdot ||v|| - (1 - t\phi(||v||)) \cdot  C_1 \cdot ||v||^2
\\
& \ge &
(C_2 - (1 - t\phi(||v||)) \cdot C_1 \cdot ||v||) \cdot ||v||
\\
& \ge &
C_2 \cdot ||v||/2.
\end{eqnarray*}

In particular we see that the only fixed point of $f$  on 
$N_{\min\{C_2/2C_1,\epsilon_2/2\},x}$ is $x$. After possibly decreasing $\epsilon_1$ 
we can assume without loss of generality
that  $N_{\epsilon_1,x} \cap N_{\epsilon_1,y} = \emptyset$ for $x,y \in
\Fix(f), x \not= y$.  Since $M$ is cocompact,
$G\backslash\Fix(f)$ is finite. 

Since $h_t = h_0$ has no fixed points outside
$D_{\min\{C_2/2C_1,\epsilon_2/2\}}T_xM$, each map $h_t$ has only one
fixed point, namely $0$. Since the $G_x$-homotopy 
$h$ is stationary outside $D_{\min\{C_2/2C_1,\epsilon_2/2\}}T_xM$,
it extends to a $G$-homotopy from $f$ to a $G$-map $f'$ such that 
$\Fix(f) = \Fix(f')$
and 
$$\exp_{x,\epsilon_1}^{-1} \circ f' \circ \exp_{x,\epsilon_2}(v) =
T_xf'(v) = T_xf(v)$$ 
holds for each $x \in \Fix (f)$ and each $v \in D_{\epsilon_3}T_xM$.
In the sequel we will identify $D_{\epsilon_1}T_xM $ with the compact 
neighborhood $N_{\epsilon_1,x}$ of $x$ by $\exp_{x,\epsilon_1}$ for $x
\in \Fix(f)$.  Since $L^{\qq G}(f)$ depends only on the $G$-homotopy class of $f$, we
can assume in the sequel 
that $f$ agrees with $T_xf\colon D_{\epsilon_3} T_xM \to
D_{\epsilon_1}T_xM$ on 
$D_{\epsilon_3} T_xM$ for each $x \in \Fix(f)$. 

Next we analyze the $G_x$-linear map $T_xf\colon T_xM \to T_x M$ for
$x \in \Fix(f)$. We can decompose the  
orthogonal $G_x$-representation $T_xM$ as 
$$T_xM ~ = ~ \bigoplus_{i=1}^n V_i^{m_i}$$
for pairwise non-isomorphic irreducible $G_x$-representations $V_1$, $V_2$,
$\ldots$ , $V_n$ and positive integers $m_1$, $m_2$, $\ldots$,  $m_n$.
The $G_x$-linear automorphism $T_xf$ splits as 
$\oplus_{i=1}^n f_i$ for $G_x$-linear automorphisms  $f_i\colon
V_i^{m_i} \to V_i^{m_i}$.  
Let $D_i = \Endo_{\rr G_x}(V_i)$ be the skew-field of $G_x$-linear
endomorphisms  
of $V_i$. It is either the field of  real numbers  $\rr$, the field of
complex numbers $\cc$ or  
the skew-field of quaternions $\hh$. There is a canonical isomorphism
of normed vector spaces
$$\Endo_{\rr G_i}(V_i^{m_i}) \cong M_{m_i}(D_i).$$
Since the open subspace $GL_{m_i}(D_i) \subseteq M_{m_i}(D_i)$ is
connected for 
$D_i = \cc$ and $D_i = \hh$ and the sign of the determinant induces a
bijection 
$\pi_0(GL_{m_i}(\rr)) \to \{\pm 1\}$, we can connect $f_i \in
\aut_{G_x}(V_i)$ 
by a (continuous) path to either $\id\colon V_i^{m_i} \to V_i^{m_i}$ or to 
$-\id_{V_i}  \oplus \id_{V_i^{m_i-1}}\colon V_i^{m_i} \to
V_i^{m_i}$. This implies 
that we can find a decomposition 
$$T_xM = V_x \oplus W_x$$
of the orthogonal $G_x$-representation $T_xM$ into orthogonal
$G_x$-subrepresentations and 
a (continuous) path $w_t\colon T_xM \to T_xM$ of linear $G_x$-maps
from $T_xf$ to 
$2 \cdot \id_{V_x} \oplus ~ 0_{W_x}$ such that $\id - w_t$ is an
isomorphism for all 
$t \in [0,1]$. Since $w_t$ is continuous on the compact set
$[0,1]$, there is a constant $C_3 \ge  1$
such that for each $v \in T_xM$ and each $t \in [0,1]$,
\begin{eqnarray*}
||w_t(v)|| & \le & C_3 \cdot ||v||.
\end{eqnarray*}
Choose a smooth function $\psi\colon [0,\epsilon_3/C_3] \to [0,1]$ such that
$\psi(t) = 1$ for $t \le \epsilon_3/3C_3$ and
$\psi(t) = 0$ for $t \ge 2\epsilon_3/3C_3$. Define a $G_x$-homotopy
$$u \colon D_{\epsilon_3/C_3}T_xM \times [0,1]  ~ \to ~ D_{\epsilon_3}T_xM,
\hspace{5mm} (y,t) ~ \mapsto ~ w_{t\cdot \psi(||v||)}(v).$$
This is a $G_x$-homotopy from $f|_{D_{\epsilon_3}T_xM} = 
T_xf|_{D_{\epsilon_3}T_xM} = w_0|_{D_{\epsilon_3}T_xM}$ to
a linear $G_x$-map $u_1$.  The map $u_1$ and the map $2\id_{V_x}
\oplus ~ 0_{W_x}\colon 
T_xM\to T_xM$ agree on $D_{\epsilon_3/3C_3}T_xM$.
For each  $t \in [0,1]$ the map
$u_t\colon D_{\epsilon_3/C_3}T_xM \to D_{\epsilon_3}T_xM$ has only one
fixed point, namely $0$, 
since this is true for $w_t$ for each $t \in [0,1]$ by construction. 
The $G_x$-homotopy $u$ is stationary outside $D_{2\epsilon/3C_3}T_xM$.
Hence it can be extended to a $G$-homotopy $U\colon M \times [0,1] \to M$
which is stationary outside $D_{2\epsilon/3C_3}T_xM$. Since 
\begin{eqnarray*}
L^{\qq G}(f)  & = & L^{\qq G}(U_1);
\\
\deg((\id - T_xf)^c) & = & \deg((\id - T_xU_1)^c),
\end{eqnarray*}
we can assume without loss of generality
that $f$ looks on $D_{\epsilon_3/C_3}T_xM$ like 
$$2\id_{V_x} \oplus ~ 0_{W_x}\colon T_xM = V_x \oplus W_x \to T_xM =
V_x \oplus W_x$$ 
for each $x \in \Fix(f)$.
By scaling the metric with a constant, we can arrange that we can take
$\epsilon/3C_3 = 1/2$ and $\epsilon_1 = 2$, in other words, 
we can identify $D_2T_xM$ with a neighborhood of $x$ in $M$ 
and $f$ is given on $D_{1/2}T_xM$ by $T_xf = 2 \cdot \id_{V_x} \oplus
0_{W_x}$.  

Let $d$ be the metric on $M$ coming from the Riemannian metric.
Choose an integer $\delta > 0$ such that 
the inequality $d(y,f(y)) \ge \delta $ holds 
for each $y \in M$, which does not lie in $D_{1/2}T_xM$ for each  $x
\in \Fix(f)$. 
Consider $x \in \Fix(f)$. Choose $G_x$-equivariant triangulations on the unit
spheres $SV_x$ and $SW_x$ such that the diameter of each simplex
measured with respect to 
the metric $d$ is smaller than $\delta/8$. Equip $[0,1]$ with the
$CW$-structure 
whose $0$-skeleton is $\{\frac{i}{2n} \mid i = 0,1,2, \ldots, 2n\}$ for some 
positive integer $n$ which will be specified later. Equip $D_1V_x$ with the
$G_x$-$CW$-structure which is induced from the product $G_x$-$CW$-structure on
$SV_x \times [0,1]$ by the quotient map 
$$SV_x \times [0,1] \to DV_x, \hspace{1mm} (y,t) \mapsto t \cdot y.$$
This is not yet the structure of a simplicial $G_x$-complex since the cells
look like cones over simplices or products of simplices. The cones over
simplices 
are again simplices and will not be changed. There is a standard
way of subdividing a product of simplices to get a simplicial structure again.
We use the resulting simplicial $G_x$-structure on $D_1V_x$. It is 
actually a $G_x$-equivariant triangulation. Define analogously a
$G_x$-simplicial structure 
on $D_1W$. 

Notice that $D_{1/2}V_x \subseteq D_1V_x$ inherits a
$G_x$-$CW$-simplicial substructure. 
We will also use a second $G_x$-simplicial structure on $D_{1/2}V_x$,
which will be denoted by $D_{1/2}V_x'$. It is induced by the product
$G_x$-$CW$-structure on 
$SV_x \times [0,1]$ above together with the the quotient map 
$$SV_x \times [0,1] \to D_{1/2} V_x, \hspace{1mm} (y,t) \mapsto t/2 \cdot y.$$
The $G_x$-simplicial -structure on $D_{1/2}V_x'$ is finer than 
the one on $D_{1/2}V_x$ but agrees with the one
on $D_{1/2}V_x$ on the boundary. The map $2\id\colon V_x \to V_x$  
induces an isomorphism of $G_x$-simplicial complexes 
$2 \id \colon D_{1/2}V_x' \xrightarrow{\cong} D_{1/2}V_x$, but it does
not induce a simplicial 
map $2\id\colon D_{1/2}V_x \to D_1V_x$. 
The latter map is at least cellular with respect to the
$G_x$-$CW$-structures induced  
from the $G_x$-simplicial structures since the $p$-skeleton of
$D_{1/2}V_x$ is contained in the $p$-skeleton of $D_{1/2}V_x'$.

We equip $D_1V_x \times D_1W_x$, $D_{1/2} V_x \times D_1W_x$ and
$D_{1/2} V_x \times D_1W_x$  with the product $G_x$-simplicial
structure. Again this requires subdividing
products of simplices (except for products of a simplex with a vertex).

Recall that we have identified $D_2T_xM$ with its image under the
exponential map. 
Choose a complete set of representatives $\{x_1, x_2, \ldots ,x_k\}$
for the $G$-orbits in $\Fix(f)$. By the construction above we get a
$G$-triangulation  on the $G$-submanifold
$\coprod_{i = 1}^k G \times_{G_{x_i}} D_1V_{x_i} \times D_1W_{x_i}$ of $M$
such that the diameter of each simplex is smaller than $\delta/4$ if
we choose the 
integer $n$ above small enough. It can be extended to a $G$-triangulation 
$K$ of $M$ such that each simplex has a diameter less than $\delta/4$.   
Let $K'$ be the refinement of $K$ which agrees with $K$ outside
$D_{1/2}V_{x_i} \times D_1W$ and is $D_{1/2}V_{x_i}' \times D_1W$ on
the subspace 
$D_{1/2}V_{x_i} \times D_1W$. Then $f\colon  K' \to K$ is a $G$-map
which is simplicial on 
$\coprod_{i = 1}^k G \times_{G_{x_i}} (D_{1/2}V_{x_i}' \times
D_1W_{x_i}T_xM)$. 
The construction in the proof of the (non-equivariant) simplicial
approximation theorem 
yields a subdivision $K''$ of $K'$ such that $K''$ and $K'$ agree on 
$\coprod_{i = 1}^k G \times_{G_{x_i}} (D_{1/2}V_{x_i}' \times
D_1W_{x_i})$ and a $G$-homotopy 
$h\colon  M \times [0,1] \to M$ from $h_0 = f$ to a simplicial map
$h_1\colon  K'' \to K$ 
such that $h$ is stationary on 
$\coprod_{i = 1}^k G \times_{G_{x_i}} (D_{1/2}V_{x_i}' \times
D_1W_{x_i})$ and  
the track of the homotopy for each point in $M$ lies within a simplex of $K$.
Recall that any simplex of $K$ has diameter at most $\delta/4$ and
$d(y,f(y)) \ge \delta$ 
holds for $y \in M$ which does not lie in 
$\coprod_{i = 1}^k G \times_{G_{x_i}} (D_{1/2}V_{x_i}' \times D_1W_{x_i})$.
Hence for any simplex $e \in K''$ outside
$\coprod_{i = 1}^k G \times_{G_{x_i}} (D_{1/2}V_{x_i}' \times
D_1W_{x_i})$ we have 
$h_1(e) \cap e = \emptyset$. The $G$-map $h_1\colon  K'' \to K''$ is
not simplicial anymore 
but at least cellular with respect to the $G$-$CW$-structure on $M$
coming from $K''$. 
This comes from the fact that  each skeleton of $K''$ is larger than
the one of $K'$. 

Next we compute $\inc(h_1,e)$ for cells $e$ in $M$ with respect to
the $G$-$CW$-structure induced by $K''$. Obviously 
$\inc(h_1,e)= 0$ if $e$ does not belong to 
$\coprod_{i = 1}^k G \times_{G_{x_i}} (D_{1/2}V_{x_i}' \times D_1W_{x_i})$ 
since for such cells $e$ we have
$h_1(e) \cap e = \emptyset$. 
If $e$ belongs to  $D_{1/2}V_{x_i}' \times D_1W_{x_i}$ its image under
$h_1 = f = 2\id_{V_{x_i}} \oplus ~ 0_{W_{x_i}}$ does not meet the interior of
$e$ unless it is the zero simplex  sitting at $(0,0) \in V_{x_i}
\oplus W_{x_i}$ or a  simplex
of the shape $\{t \cdot x\mid t \in [0,1/4n], x \in e \} \times \{0\}$
for some simplex in $e \in 
SV_{x_i}$. Hence among the cells in $D_{1/2}V_{x_i}' \times D_1W_{x_i}$
only the zero simplex  sitting at $(0,0) \in V_{x_i} \oplus W_{x_i}$
and the simplex 
of the shape $\{t \cdot x\mid t \in [0,1/4n], x \in e \} \times \{0\}$
for some simplex in $e \in 
SV_{x_i}$ can have non-zero incidence numbers $\inc(f,e)$ and one
easily checks that these incidence numbers are all equal to $1$.
Hence we get using Lemma \ref{lem: decription of L^{QG}(f) in terms of incidence numbers}
and the equality $\inc(f,e) = \inc(f,ge)$
\begin{eqnarray*}
L^{\qq G}(f) & = & L^{\qq G}(h_1)
\\
& = & 
\sum_{p \ge 0} (-1)^p \cdot \sum_{Ge \in G\backslash I_p(K'')} |G_e|^{-1} \cdot \inc(f,e)
\\
& = & 
\sum_{p \ge 0} (-1)^p \cdot \sum_{i=1}^k ~ 
\sum_{\substack{G_{x_i}e \in \\G_{x_i}\backslash I_p(D_{1/2}V_{x_i}' \times D_1W)}} 
|(G_{x_i})_e|^{-1} \cdot \inc(f,e)
\\ 
& = & 
\sum_{p \ge 0} (-1)^p \cdot \sum_{i=1}^k |G_{x_i}|^{-1} 
\sum_{\substack{G_{x_i}e \in \\G_{x_i}\backslash I_p(D_{1/2}V_{x_i}' \times D_1W)}} 
|G_{x_i}/(G_{x_i})_e| \cdot \inc(f,e)
\\
& = &  
\sum_{p \ge 0} (-1)^p \cdot \sum_{i=1}^k |G_{x_i}|^{-1} 
\sum_{e \in I_p(D_{1/2}V_{x_i}' \times D_1W)} 
\inc(f,e)
\\ 
& = & 
\sum_{i=1}^k |G_{x_i}|^{-1} \sum_{p \ge 0} (-1)^p \cdot 
\sum_{e \in I_p(D_{1/2}V_{x_i}' \times D_1W)}  \inc(f,e)
\\ 
& = & 
\sum_{i=1}^k |G_{x_i}|^{-1} \left( 1 + \sum_{p \ge 1} (-1)^p \cdot 
    |I_{p-1}(SV_{x_i})| \right)
\\
& = & 
\sum_{i=1}^k |G_{x_i}|^{-1} \left( 1 -  \chi(SV_{x_i})\right)
\\
& = & 
\sum_{i=1}^k |G_{x_i}|^{-1} (-1)^{\dim(V_{x_i})} 
\\
& = & 
\sum_{i=1}^k |G_{x_i}|^{-1} \frac{\det(\id - T_{x_i}f)}{|\det(\id - T_{x_i}f)|}
\\
& = &
\sum_{i=1}^k |G_{x_i}|^{-1} \deg((\id - T_{x_i}f)^c)
\\
& = & 
\sum_{Gx \in G\backslash \Fix(f)} |G_x|^{-1} \cdot \deg((\id - T_{x_i}f)^c).
\end{eqnarray*}

This finishes the proof of Theorem
\ref{the: the orbifold Lefschetz fixed point theorem}. \qed

%%%%%%%%%%%%%%%%%%%%%%%%%%   Section 3 %%%%%%%%%%%%%%%%%%%%%%%%%% %%%%%%%%%%%%

\typeout{-----------------------  Section 3  ------------------------}

\tit{The equivariant Lefschetz classes}{The equivariant Lefschetz classes}

In this section we define the equivariant Lefschetz class appearing
in the equivariant Lefschetz fixed point Theorem
\ref{the: equivariant Lefschetz theorem}.
We will use the following notation in the sequel.

\begin{notation} \label{not: X^H(x) and so on}
Let $G$ be a discrete group and $H \subseteq G$ be a subgroup.
Let $N\!H = \{g \in G \mid gHg^{-1} = H\}$ be its \emph{normalizer}
and let $W\!H := N\!H/H$ be its \emph{Weyl group}.

Denote by $\consub(G)$ the set of conjugacy classes $(H)$ of subgroups
$H \subseteq G$. 

Let $X$ be a $G$-$CW$-complex. Put 
\begin{eqnarray*}
X^H & := & \{x \in X \mid H \subseteq G_x\};
\\
X^{>H} & := & \{x \in X \mid H \subsetneq G_x\},
\end{eqnarray*} 
where $G_x$ is the isotropy group of $x$ under the $G$-action.

Let $x\colon G/H \to X$ be a $G$-map. Let $X^H(x)$ be the component of
$X^H$ containing $x(1H)$. Put 
$$X^{>H}(x) = X^H(x) \cap X^{>H}.$$
Let $W\!H_x$ be the isotropy group of $X^H(x) \in \pi_0(X^H)$ under the $W\!H$-action.
\end{notation}

Next we define the group $U^G(X)$, where the equivariant Lefschetz
class will take its values.

Let $\Pi_0(G,X)$ be the \emph{component category} of the $G$-space $X$
in the sense of  tom Dieck \cite[I.10.3]{Dieck(1987)}.
Objects are $G$-maps $x\colon G/H \to X$. A morphism $\sigma$ from
$x\colon G/H \to X$ to 
$y\colon G/K \to X$ is a $G$-map $\sigma\colon G/H \to G/K$ 
such that $y \circ \sigma$ and $x$ are $G$-homotopic. A $G$-map
$f\colon  X \to Y$ induces a functor $\Pi_0(G,f)\colon  \Pi_0(G,X) \to
\Pi_0(G,Y)$ 
by composition with $f$. Denote by $\Is \Pi_0(G,X)$ the set of
isomorphism classes $[x]$ of objects $x\colon G/H \to X$ in
$\Pi_0(G,X)$. Define  
\begin{eqnarray}
U^G(X) & := & \zz[\Is \Pi_0(G,X)], 
\label{definition of U^G(X)}
\end{eqnarray}
where for a set $S$ we denote by $\zz [S]$ 
the free abelian group with basis $S$.
Thus we obtain a covariant functor from the category of $G$-spaces to the 
category of abelian groups. Obviously $U^G(f) = U^G(g)$ if $f,g \colon
X \to Y$ are $G$-homotopic. 

There is a natural bijection
\begin{eqnarray}
\Is \Pi_0(G,X)  & \xrightarrow{\cong} & \coprod_{(H) \in \consub(G)} 
W\!H\backslash \pi_0(X^H),
\label{identifying Is Pi_0(G,X)}
\end{eqnarray}
which sends $x\colon G/H \to X$ to the orbit under the $W\!H$-action on $\pi_0(X^H)$
of the component $X^H(x)$ of $X^H$ which contains the point $x(1H)$.
It induces a natural isomorphism
\begin{eqnarray}
U^G(X)  & \xrightarrow{\cong} & \bigoplus_{(H) \in \consub(G)} 
\zz [W\!H\backslash \pi_0(X^H)]. 
\label{identifying U^G(X)}
\end{eqnarray}

Let $\alpha \colon G \to K $ be a group homomorphism and $X$
be a $G$-$CW$-complex. We obtain from $\alpha$ a functor
$$\alpha_*\colon \Pi_0(G,X) \to \Pi_0(K,\alpha_* X)$$
which sends an object $x\colon G/H \to X$ to the object
$\alpha_*(x)\colon  K/\alpha (H) = \alpha_*(G/H) \to \alpha_*X$ and
similarly for morphisms. Thus 
we obtain an \emph{induction homomorphism} of abelian groups
\begin{eqnarray}
\alpha_*\colon U^G(X) & \to & U^K(\alpha_*X).
\label{alpha_*: U^G(X) to U^K(alpha_*X)}
\end{eqnarray}

Next we define the equivariant Lefschetz class.
Let $X$ be a finite proper $G$-$CW$-complex.
Let $f\colon X \to X$ be a cellular $G$-map such that for each subgroup $K
\subseteq G$ the map  
$\pi_0(f^K)\colon \pi_0(X^K) \to \pi_0(X^K)$ is the identity.  
For any $G$-map $x\colon  G/H \to X$ it induces a map 
$$(f^H(x),f^{>H}(x)) \colon (X^H(x),X^{>H}(x)) \to (X^H(x),X^{>H}(x))$$
of pairs of finite proper $W\!H_x$-$CW$-complexes. Then 
$$L^{\zz W\!H_x}(f^H(x),f^{>H}(x)) \in \zz$$
is defined (see \eqref{definition of L^{RG}(f)})
since the isotropy group under the $W\!H_x$-action of any
point in $X^H(x) - X^{>H}(x)$ is trivial. 

\begin{definition} \label{def: equivariant Lefschetz class Lambda^G(f)}
We define the 
\emph{equivariant Lefschetz class} of $f$
\begin{eqnarray*}
\Lambda^G(f) & \in & U^G(X)
\end{eqnarray*}
by assigning  to $[x\colon G/H \to X] \in \Is \Pi_0(G,X)$ the integer 
$$L^{\zz W\!H_x}\left(f^H(x),f^{>H}(x)\colon (X^H(x), X^{>H}(x)) \to (X^H(x), X^{>H}(x))\right),$$ 
if $f^H\colon X^H \to X^H$ maps $X^H(x)$ to itself, and zero otherwise.
\end{definition}

Since $X^{>H}(x) \not= X^H(x)$ and therefore $L^{\zz W\!H_x}(f^H(x),f^{>H}(x)) \not= 0$
holds only for finitely many elements $[x]$ in $\Is \Pi_0(G,X)$,
Definition \ref{def: equivariant Lefschetz class Lambda^G(f)} makes sense.
Notice for the sequel that $f^H(X^H(x)) \cap X^H(x) \not= \emptyset$ implies
$f^H(X^H(x)) \subseteq X^H(x)$. The elementary proof that 
Lemma \ref{lem: properties of L^{RG}(f)} implies the following lemma 
is left to the reader. 

\begin{lemma} 
\label{lem: properties of Lambda^{G}(f)}
Let $X$ be a finite proper $G$-$CW$-complex.
Let $f\colon X \to X$ be a cellular $G$-map. Then

\begin{enumerate}
\item \label{lem: properties of Lambda^{G}(f): G-homotopy invariance}
The equivariant Lefschetz class $\Lambda^G(f)$ depends only on the
cellular $G$-homotopy class of $f$; 

\item \label{lem: properties of Lambda^{G}(f): conjugation invariance}
If $f'\colon  Y \to Y$ is a cellular $G$-self-map of a finite $G$-$CW$-complex
and $h\colon  X \to Y$ is a cellular $G$-homotopy equivalence satisfying
$h \circ f \simeq_G f' \circ h$, then $U^G(h)\colon  U^G(X)
\xrightarrow{\cong}  U^G(Y)$  
is bijective and sends $\Lambda^G(f)$ to $\Lambda^G(f')$;

\item \label{lem: properties of Lambda^{G}(f): induction}
Let $\alpha\colon  G \to K$ be an inclusion of groups. Denote by
$\alpha_*f$ the cellular $K$-self-map obtained by induction with $\alpha$. Then
$$\Lambda^{K}(\alpha_*f) ~ = ~ \alpha_*\Lambda^{G}(f).$$

\end{enumerate}
\end{lemma}

By the equivariant cellular approximation theorem (see for instance
\cite[Theorem 2.1 on page 32]{Lueck(1989)}) any $G$-map of
$G$-$CW$-complexes is 
$G$-homotopic to a cellular $G$-map and two cellular $G$-maps which
are $G$-homotopic are actually cellularly $G$-homotopic. 
Hence we can drop the assumption cellular
in the sequel because of $G$-homotopy invariance of the equivariant
Lefschetz class 
(see Lemma \ref{lem: properties of Lambda^{G}(f)}
\eqref{lem: properties of Lambda^{G}(f): G-homotopy invariance}).

%%%%%%%%%%%%%%%%%%%%%%%%%%   Section 4 %%%%%%%%%%%%%%%%%%%%%%%%%%%%%%%%%%%%%%%%

\typeout{-----------------------  Section 4  ------------------------}

\tit{The local equivariant Lefschetz class}{The local equivariant
Lefschetz class} 

In this section we introduce the local equivariant Lefschetz class in
terms of fixed point data. 
Before we can define it, we recall the classical notion of the 
equivariant Lefschetz class with values in the Burnside ring for a
finite group. 

Let $K$ be a finite group. The abelian group $U^K(\{\ast\})$ is
canonically isomorphic to the abelian group which underlies the
Burnside ring $A(K)$. Recall that the \emph{Burnside ring} is the
Grothendieck ring of finite $K$-sets with the additive structure
coming from disjoint union and the multiplicative structure coming
from the Cartesian product. 

Let $X$ be a finite $K$-$CW$-complex. Define the 
\emph{equivariant Lefschetz class with values in the Burnside ring} of $f$
\begin{eqnarray}
\Lambda^K_0(f) & \in & A(K) = U^K(\{\ast\})
\label{definition of Lambda^K_0(f) for finite K}
\end{eqnarray} 
by
$$\Lambda^K_0(f) ~ := ~ \sum_{(H) \in \consub(K)}
L^{\zz W\!H}(f^H,f^{>H}) \cdot [K/H].$$
(Here and elsewhere the subscript $_0$ shall indicate that
the corresponding invariant takes values in the Burnside ring 
and the component structure of the various fixed point sets is not taken
into account.)
Denote by 
\begin{eqnarray}
\ch^K_0 \colon  A(K) & \to & \prod_{(H) \in \consub(K)} \zz
\label{character map ch^K_0 of A(K) for finite K}
\end{eqnarray}
the \emph{character map} which sends the class of a finite set $S$ to the
collection $\{|S^H|\mid (H) \in \consub(K)\}$ given by the orders of
the various $H$-fixed point sets. The character map is a
ring homomorphism, and it is injective (see 
Lemma \ref{lem: injectivity of the character map chi^G}). The
equivariant Lefschetz class  
$\Lambda^K_0(f)$ is characterized by the property
(see for instance 
\cite[Theorem 2.19 on page 504]{Laitinen-Lueck(1989)}), 
\cite[Lemma 3.3 on page 138]{Lueck(1988)})
\begin{eqnarray}
\ch^K_0(\Lambda_0^K(f)) & = & \{L^{\zz[\{1\}]}(f^H) \mid (H) \in \consub(K)\}.
\label{Lambda^G_0 and the character map}
\end{eqnarray}
If $\pr\colon X \to \{\ast\}$ is the projection, then
$U^K(\pr)\colon U^K(X) \to U^K(\{\ast\}) = A(K)$ sends $\Lambda^K(f)$
(see Definition \ref{def: equivariant Lefschetz class Lambda^G(f)}) to
$\Lambda_0^K(f)$ defined in \eqref{definition of Lambda^K_0(f) for finite K}.

Let $V$ be a (finite-dimensional) $K$-representation and let 
$f\colon V^c \to V^c$ be a $K$-self-map of the
one-point-compactification $V^c$. Define its \emph{equivariant degree}
\begin{eqnarray}
\Deg_0^K(f) & \in  & A(K) = U^{K}(\{\ast\})
\label{definition of Deg^K(f)}
\end{eqnarray}
by
$$\Deg_0^K(f) ~ := ~ (\Lambda_0^K(f) - 1) \cdot
(\Lambda_0^K(\id_{V^c}) - 1).$$
Since the character map \eqref{character map ch^K_0 of A(K) for finite K} is
an injective ring homomorphism, we conclude from
 \eqref{Lambda^G_0 and the character map} above that $\Deg_0^K(f)$  is
uniquely characterized by the equality
\begin{eqnarray}
\ch^K_0(\Deg_0^K(f)) & = & \{\deg(f^H) \mid (H) \in \consub(K)\},
\label{Lambda^G and the character map}
\end{eqnarray}
where $\deg(f^H)$ is the degree of the self-map 
$f^H\colon (V^c)^H \to (V^c)^H$ of the connected closed orientable
manifold $(V^c)^H$, 
if $\dim((V^c)^H) \ge 1$, and $\deg(f^H)$ is defined to be $1$, 
if $\dim((V^c)^H) = 0$. The equivariant degree of \eqref{definition of
Deg^K(f)} induces 
an isomorphism from the $K$-equivariant stable cohomotopy of a point to
the Burnside ring $A(K)$ \cite[Theorem 7.6.7 on page 190]{Dieck(1979)},
\cite{Segal(1970)}.

Let $M$ be a cocompact
proper $G$-manifold (possibly with boundary). Let 
$f\colon  M \to M$ be a smooth $G$-map.
Denote by $\Fix(f) = \{x \in X \mid f(x) = x\}$ the 
\emph{set of fixed points} of $f$. 
Suppose that for any $x \in \Fix(f)$ the determinant
of the linear map $\id - T_xf\colon  T_xM \to T_xM$ is different from
zero. (One can always find a representative in the $G$-homotopy class
of $f$ which satisfies this assumption.) Then $G\backslash\Fix(f)$ is finite.
Consider an element $x \in \Fix(f)$.
Let $\alpha_x\colon G_x \to G$ be the inclusion. We obtain from
$(\alpha_x)_*$ (see \eqref{alpha_*: U^G(X) to U^K(alpha_*X)}) and
$U^G(x)$ for $x \in \Fix(f)$ interpreted as a $G$-map $x\colon  G/G_x
\to X$ a homomorphism 
$$U^{G_x}(\{\ast\}) \xrightarrow{(\alpha_x)_*} U^G(G/G_x)
\xrightarrow{U^G(x)} U^G(X).$$ 
Thus we can assign to $x \in \Fix(f) $ the element
$U^G(x) \circ (\alpha_x)_*(\Deg_0^{G_x}((\id - T_xf)^c))$, where 
$\Deg_0^{G_x}((\id - T_xf)^c)$ is the equivariant degree (see
\eqref{definition of Deg^K(f)}) of the map induced on the
one-point-compactifications 
by the isomorphism $(\id - T_xf)\colon T_xM \linebreak \to T_xM$.
One easily checks that this element depends only on the $G$-orbit of $x
\in \Fix(X)$. 

\begin{definition} \label{def: Lambda^G_{loc}(f)}
We can define the \emph{local equivariant Lefschetz class} by
\begin{multline*}
\Lambda_{\loc}^G(f)
\\ := ~
\sum_{Gx \in G\backslash\Fix(f)}~ 
U^G(x) \circ (\alpha_x)_*\left(\Deg_0^{G_x}\left((\id - T_xf)^c)\right)\right)
\hspace{2mm} \in U^G(M).
\end{multline*}
\end{definition}
Now have defined all the ingredients appearing in the 
Equivariant Lefschetz fixed point theorem
\ref{the: equivariant Lefschetz theorem}.
Before we give its proof, we discuss the following example

\begin{example} \label{exa: Lefschetz fixed point theorem for all
isotropy groups of odd order} 
\em 
Let $G$ be a discrete group and $M$ be a cocompact proper $G$-manifold
(possibly with boundary). 
Suppose that the isotropy group $G_x$ of each point $x \in M$ has odd order.
This holds  automatically  if $G$ itself is a finite group of odd order.
Let $f\colon  M \to  M$ be a smooth $G$-map. 
Suppose that $\Fix(f) \cap \partial M = \emptyset$ and for
each $x \in \Fix(f)$ the determinant 
of the linear map $\id - T_xf\colon  T_xM \to T_xM$ is different from
zero. If $H$ is a finite group of odd order, then the multiplicative
group of units 
$A(H)^*$ of the Burnside ring is known to be $\{\pm 1\}$
\cite[Proposition 1.5.1]{Dieck(1979)}.
The element $\Deg_0^{G_x}((\id - T_xf)^c) \in A(G_x) =
U^{G_x}(\{\ast\})$ satisfies 
$\left(\Deg_0^{G_x}((\id - T_xf)^c)\right)^2 = 1$ since this holds for
its image under the injective 
ring homomorphism $\ch^{G_x}_0\colon A(G_x) \to \prod_{(H) 
\in \consub(G_x)} \zz$,
whose coefficient at $(H) \in \consub(G_x)$ is 
$\deg((\id - T_xf)^c)^H) \in \{\pm 1\}$ (see \eqref{Lambda^G and the
character map}). Hence 
$\Deg_0^{G_x}((\id - T_xf)^c)$ belongs to  $A(G_x)^* = \{\pm
1\}$. This implies that 
$$\Deg_0^{G_x}((\id - T_xf)^c) = 
\frac{\det(\id - T_xf\colon  T_xM \to T_xM)}{|\det(\id - T_xf\colon
T_xM \to T_xM)|}  
\cdot [G_x/G_x].$$
Hence the definition of the local equivariant Lefschetz class reduces to
$$
\Lambda_{\loc}^G(f)
~ := ~
\sum_{Gx \in G\backslash\Fix(f)}~ 
\frac{\det(\id - T_xf\colon  T_xM \to T_xM)}{|\det(\id - T_xf\colon  T_xM \to T_xM)|} 
\cdot [x\colon G/G_x \to M].
$$
where $x\colon G/G_x \to M$ sends $g\cdot G_x$ to $gx$. 
\em
\end{example}

\begin{remark} \label{rem: further lefschetz classes} \em
Equivariant Lefschetz classes for compact Lie groups were studied in
\cite{Laitinen-Lueck(1989)}. In the non-equivariant setting, 
universal Lefschetz classes with values in certain $K$-groups were
defined and analyzed in \cite{Lueck(1999)}. 
It seems to be possible to combine the $K$-theoretic invariants there
with the equivariant versions presented here to obtain
a universal equivariant Lefschetz class.
\em
\end{remark}

%%%%%%%%%%%%%%%%%%%%%%%%%%   Section 5 %%%%%%%%%%%%%%%%%%%%%%%%%%%%%%%%%%%%%%%%

\typeout{-----------------------  Section 5  ------------------------}

\tit{The proof of the equivariant Lefschetz fixed point theorem}
{The proof of the equivariant Lefschetz fixed point theorem}

This section is devoted to the proof of the 
equivariant Lefschetz fixed point Theorem 
\ref{the: equivariant Lefschetz theorem}.

First we define the \emph{character map} for a proper $G$-$CW$-complex $X$:
\begin{eqnarray}
\ch^G(X)\colon U^G(X) & \to & \bigoplus_{\Is \Pi_0(G,X)}\qq.
\label{definition of character map chi^G(X)}
\end{eqnarray}

We have to define for an isomorphism class $[x]$ of objects $x\colon
G/H \to X$ in  
$\Pi_0(G,X)$ the component $\ch^G(X)([x])_{[y]}$ of $\ch^G(X)([x])$
which belongs to an isomorphism class $[y]$ of objects $y\colon G/K \to X$ in 
$\Pi_0(G,X)$, and check that $\chi^G(X)([x])_{[y]}$ is different from zero for
at most finitely many $[y]$. Denote by $\mor(y,x)$ the set of morphisms from $y$ to $x$ in
$\Pi_0(G,X)$. We have the left operation 
$$\aut(y) \times \mor(y,x) \to \mor(y,x), 
\hspace{5mm}
(\sigma, \tau) \mapsto \tau \circ \sigma^{-1}.$$
There is an isomorphism of groups
$$W\!K_y \xrightarrow{\cong} \aut(y)$$
which sends $gK \in W\!K_y$ to  the automorphism of $y$ given by the $G$-map
$$R_{g^{-1}}\colon G/K \to G/K,
\hspace{5mm}
g'K \mapsto g'g^{-1}K.$$
Thus $\mor(y,x)$ becomes a left $W\!K_y$-set. 

The $W\!K_y$-set $\mor(y,x)$ can be rewritten as
$$\mor(y,x) ~ = ~ \{g \in G/H^K \mid g\cdot x(1H) \in X^K(y)\},$$
where the left operation of $W\!K_y$ on 
$\{g \in G/H^K \mid g\cdot x(1H) \in Y^K(y)\}$ comes from the canonical
left action of $G$ on $G/H$. Since $H$ is finite and hence
contains only finitely many subgroups, the set $W\!K\backslash
(G/H^K)$ is finite 
for each $K \subseteq G$ and is non-empty for only finitely many
conjugacy classes $(K)$ 
of subgroups $K \subseteq G$.
This shows that $\mor(y,x) \not= \emptyset$  for at most finitely many isomorphism classes $[y]$
of objects $y \in \Pi_0(G,X)$ and that the $W\!K_y$-set $\mor(y,x)$ decomposes into finitely many 
$W\!K_y$ orbits with finite isotropy groups for each object $y \in \Pi_0(G,X)$. We define
\begin{eqnarray}
\ch^G(X)([x])_{[y]} & := & \sum_{\substack{W\!K_y \cdot \sigma \in\\
W\!K_y\backslash \mor(y,x)}}~  |(W\!K_y)_{\sigma}|^{-1},
\label{definition of chi^G(X)([x])_{[y]}}
\end{eqnarray}
where $(W\!K_y)_{\sigma}$ is the isotropy group of $\sigma \in
\mor(y,x)$ under the $W\!K_y$-action. Notice
that $\ch^G(X)([x])_{[y]}$ is the same as
$\dim_{\qq W\!K_y}(\qq(\mor(y,x)))$, 
if one defines $\dim_{\qq W\!K_y}(P)$ of a finitely generated $\qq W\!K_y$-module $P$ 
by $\tr_{\qq W\!K_y}(\id_P)$ (see Lemma \ref{lem: basic properties of tr_{qq G}(a)}
\eqref{lem: basic properties of tr_{qq G}(a): dimension and trace}).
This agrees with the more general notion of
the von Neumann dimension of the finitely generated Hilbert $\caln(W\!K_y)$-module
$l^2(\mor(y,x))$.

The character map $\ch^G(X)$ of 
\eqref{definition of character map chi^G(X)} 
should not be confused with the isomorphism appearing in 
\eqref{identifying U^G(X)}. If $G$ is finite and $X = \{\ast\}$, then
character map $\ch^G(X)$ of 
\eqref{definition of character map chi^G(X)} and the character map
$\ch^G_0$ of \eqref{character map ch^K_0 of A(K) for finite K}  are related
under the identifications $U^G(\{\ast\}) = A(G)$  and
$\Is\Pi_0(G,\{\ast\}) = \consub(G)$ by
$$(\ch^G_0)_{(H)} ~ = ~ |W\!H| \cdot \ch^G(\{\ast\})_{(H)}$$
for $(H) \in \consub(G)$.

\begin{lemma} 
\label{lem: injectivity of the character map chi^G}
The map $\ch^G$ of \eqref{definition of character map chi^G(X)} is injective.
\end{lemma}
\proof
Consider $u \in U^G(X)$ with $\ch^G(X)(u) = 0$.
We can write $u$ as a finite sum 
$$u = \sum_{i=1}^n m_i \cdot [x_i\colon G/H_i \to X]$$
for some integer $n \ge 1$ and integers $m_i$ such that 
$[x_i] = [x_j]$ implies $i = j$ and such that 
$H_i$ is subconjugate to $H_j$ only if $i \ge j$ or $(H_i) =
(H_j)$. We have to show 
that $u = 0$. It suffices to prove 
$$\ch^G(X)([x_i])_{[x_1]} ~ = ~
\left\{
\begin{array}{lll} 0 & & \mbox{ if } i > 1
\\
1 & & \mbox{ if } i = 1
\end{array}
\right.
$$
because then $m_1 = 0$ follows and one can show inductively
that $m_i = 0$ for $i = 1,2, \ldots , n$.
Suppose that $\ch^G(X)([x_i])_{[x_1]} \not= 0$. 
Then $\mor(x_1,x_i)$ is non-empty. This implies
that $\im(x_i) \cap X^{H_1}$ is non-empty and hence that 
$H_1$ is subconjugate to $H_i$. From the way we have enumerated the
$H_i$'s we conclude $(H_i) = (H_1)$. Since $\im(x_i) \cap X^{H_1}(x_1)$
is non-empty, we get $[x_i] = [x_1]$ and hence $i = 1$. Since by definition
$\ch^G(X)([x_1])_{[x_1]} ~ = ~ 1$,
the claim follows. \qed 

\begin{lemma} \label{character map and Lefschetz number}
Let $f\colon X \to X$ be a $G$-self-map of a finite proper
$G$-$CW$-complex $X$. Let 
$[y]$ be an isomorphism class of objects $y\colon G/K \to X$ in $\Pi_0(G,X)$.
Then
$$\ch^G(X)(\Lambda^G(f))_{[y]} ~ = ~ 
L^{\qq W\!K_y}\left(f|_{X^K(y)}\colon X^K(y)  \to X^K(y)\right),$$
if $f^K(X^K(y)) \subseteq X^K(y)$ and
$$\ch^G(X)(\Lambda^G(f))_{[y]} ~ = ~ 0$$
otherwise.
\end{lemma}
\proof We first consider the case $f^K(X^K(y)) \subseteq X^K(y)$.
Let $X_p$ be the
$p$-th skeleton of $X$. Then we can write $X_p$ as a $G$-pushout
\comsquare{\coprod_{i=1}^{n_p} G/H_i \times S^{p-1}}
{\coprod_{i=1}^{n_p} q_{p,i}}{X_{p-1}}
{}{}
{\coprod_{i=1}^{n_p} G/H_i \times D^p}
{\coprod_{i=1}^{n_p} Q_{p,i}}{X_p}
for an integer $n_p \ge 0$ and finite subgroups $H_i \subseteq G$.
For each $i \in I_p$ let $x_{p,i}\colon G/H_i \to X$ be the $G$-map
obtained by 
restricting the characteristic map $Q_{p,i}$ to $G/H_i \times
\{0\}$. For $i = 1,2, \ldots, p$ define
\begin{eqnarray}
\inc(f,p,i) & = & \inc(f,e_i)
\label{def of inc(f,p,i)}
\end{eqnarray}
where $e_i$ is the open cell $Q_{p,i}(gH_i \times (D^p - S^{p-1}))$ for any
choice of $gH_i \in G/H_i$ 
and $\inc(f,e_i)$ is the incidence number defined in 
\eqref{definition of incidence number inc(f,e)}.
Since $f$ is $G$-equivariant, the choice of $gH_i \in G/H_i$ does not matter.

Now the $G$-$CW$-structure on $X$ induces a $W\!K_y$-$CW$-structure on
$X^K(y)$ whose $p$-skeleton $X^K(y)_p$ is $X^K(y)\cap X_p$. 
Let $x\colon G/H \to X$ be an object in $\Pi_0(G,X)$.
Recall that $\mor(y,x)$ is the set of morphisms from $y$ to $x$ in
$\Pi_0(G,X)$ and carries a canonical $W\!K_y$-operation.
One easily checks that
there is a $W\!K_y$-pushout diagram
\comsquare{\coprod_{i=1}^{n_p} \mor(y,x_{p,i}) \times S^{p-1}}
{\coprod_{i=1}^{n_p} q_{p,i}^K(y)}{X^K(y)_{p-1}}
{}{}
{\coprod_{i=1}^{n_p} \mor(y,x_{p,i}) \times D^p}
{\coprod_{i=1}^{n_p} Q_{p,i}^K(y)}{\,\,X^K(y)_p\,,}
where the maps $q_{p,i}^K(y)$ and $Q_{p,i}^K(y)$ are induced by the maps
$q_{p,i}$ and $Q_{p,i}$ in the obvious way. 
We conclude from Lemma \ref{lem: decription of L^{QG}(f) in terms of
incidence numbers} 
\begin{multline}
L^{\qq[W\!K_y]}(f^K(y))
\\ 
= ~ 
\sum_{p \ge 0} (-1)^p \cdot \sum_{i=1}^{n_p} ~ 
\sum_{\substack{W\!K_y \cdot \sigma \in \\
W\!K_y\backslash \mor(y,x_{p,i})}}~ 
|(W\!K_y)_{\sigma}|^{-1} \cdot \inc(f,p,i),
\label{L^{qq[WK_y]}(f^K(y)) and inc}
\end{multline}
where the incidence number $\inc(f,p,i)$ has been defined in
\eqref{def of inc(f,p,i)}.

Analogously we get for any object $x\colon G/H \to X$ a
$W\!H_x$-pushout diagram
\comsquare{\coprod_{\substack{i=1,2 \ldots, n_p\\{[x_{p,i}] = [x]}}}
\aut(x,x) \times S^{p-1}} 
{\coprod_{i=1}^{n_p} q_{p,i}^H(x)}{X^H(x)_{p-1} \cup X^{>H}(x)}
{}{}
{\coprod_{\substack{i=1,2 \ldots, n_p\\{[x_{p,i}] = [x]}}} \aut(x,x)
\times D^p} 
{\coprod_{i=1}^{n_p} Q_{p,i}^H(x)}{X^H(x)_p\cup X^{>H}(x)}
describing how the $p$-skeleton of the relative $W\!H_x$-$CW$-complex
$(X^H(x)$, $X^{>H}(x))$ is obtained from its $(p-1)$-skeleton.
We conclude from Lemma \ref{lem: decription of L^{QG}(f) in terms of
incidence numbers} 
\begin{multline}
L^{\zz[W\!H_x]}\left(f^H(x),f^{>H}(x))\colon (X^H(x),X^{>H}(x)) \to
  (X^H(x),X^{>H}(x))\right) 
\\
~ = ~ 
\sum_{p \ge 0} (-1)^p \cdot \sum_{\substack{i = 1,2, \ldots
  ,n_p\\{[x_{p,i}] = [x]}}} 
  \inc(f,p,i),
\label{L^{zz[WK_y]}(f^K(y),f^{>K}(y))}
\end{multline}
where $\inc(f,p,i)$ is the number introduced in  
\eqref{def of inc(f,p,i)}.
We get  from the Definition 
\ref{def: equivariant Lefschetz class Lambda^G(f)}
of $\Lambda^G(f)$ and \eqref{L^{zz[WK_y]}(f^K(y),f^{>K}(y))}
\begin{eqnarray}
\hspace{-5mm} \ch^G(X)(\Lambda^G(f))_{[y]} & = & 
\sum_{p \ge 0} (-1)^p \cdot \sum_{i=1}^{n_p}  \inc(f,p,i) \cdot
\ch^G(X)([x_{p,i}])_{[y]}. 
\label{chi^G(X)(Lambda^G(f))_{[y]}}
\end{eqnarray}
Now Lemma \ref{character map and Lefschetz number} follows from
the definition \eqref{definition of chi^G(X)([x])_{[y]}} of
$\ch^G(X)([x_{p,i}])_{[y]}$ and equations 
\eqref{L^{qq[WK_y]}(f^K(y)) and inc} and
\eqref{chi^G(X)(Lambda^G(f))_{[y]}} in the case $f^K(X^K(y)) \subseteq
X^K(y)$. 

Now suppose that $f^K(X^K(y)) \cap  X^K(y) = \emptyset$. For any object $x\colon G/H \to X$ with
$\mor(y,x) \not= \emptyset$ we conclude $f^H(X^H(x)) \cap  X^H(x) = \emptyset$ and hence
$\Lambda^G(f)$ assigns to $[x]$ zero by definition. We conclude 
 $\ch^G(X)(\Lambda^G(f))_{[y]} ~ = ~ 0$ from the definition of $\ch^G$. 
This finishes the proof of Lemma \ref{character map and Lefschetz number}.
\qed

\begin{lemma}
\label{character map and local Lefschetz class}
Let $M$ be a cocompact smooth proper $G$-manifold. Let $f\colon M\to
M$ be a smooth 
$G$-map. Suppose that for any $x \in Fix(f)$ the determinant
$\det(\id_{T_xf} - T_xf)$ is different from zero. Let $y\colon G/K \to
M$ be an object in 
$\Pi_0(G,M)$. Denote by $(W\!K_y)_x$ the isotropy group for
$x \in \Fix(f|_{M^K(y)})$ under the $W\!K_y$-action 
on $M^K(y)$. 

Then $G\backslash\Fix(f)$ is finite and we get
\begin{multline*}
\ch^G(M)_{[y]}(\Lambda_{\loc}^G(f)) \\
~ = ~ 
\sum_{\substack{W\!K_y\cdot x \in \\W\!K_y\backslash
\Fix(f|_{M^K(y)})}} \left|(W\!K_y)_x\right|^{-1} \cdot 
\deg\left(\left(\id_{T_xM^K(y)} - T_x(f|_{M^K(y)})\right)^c\right).
\end{multline*}
\end{lemma}
\proof
We have already shown in Theorem \ref{the: the orbifold Lefschetz fixed point theorem}
that $G\backslash\Fix(f)$ is finite.

Consider $x \in M$. Let $G_x$ be its isotropy group under the
$G$-action and denote 
by $x\colon G/G_x \to M$ the $G$-map $g \mapsto gx$. Let $\alpha_x\colon G_x \to
G$ be the inclusion. 
Given a morphism $\sigma\colon y \to x$, let $(W\!K_y)_{\sigma}$ be
its isotropy group under 
the $W\!K_y$-action and let $(K_{\sigma})$ be the element in
$\consub(G_x)$ given by 
$g^{-1}Kg$ for any element $g \in G$ for which $\sigma\colon G/K \to G/G_x$ 
sends $g'K$ to $g'gG_x$. Notice that $(K_{\sigma})$ depends only
on $W\!K_y\cdot \sigma \in W\!K_y\backslash \mor(y,x)$.
We first show for the composition
$$U^{G_x}(\{\ast\}) \xrightarrow{(\alpha_x)_*} U^G(G/G_x)
\xrightarrow{U^G(x)} U^G(M) 
\xrightarrow{\ch^G(M)_{[y]}} \qq$$
that for each $u \in U^{G_x}(\{\ast\})$
\begin{multline}
\ch^G(M)_{[y]} \circ U^G(x) \circ (\alpha_x)_* (u) \\
~ = ~ \sum_{\substack{W\!K_y\cdot \sigma \in \\W\!K_y\backslash \mor(y,x)}}~
|(W\!K_y)_{\sigma}|^{-1} \cdot \ch^{G_x}_0(u)_{(K_{\sigma})}.
\label{ch^G(M)_[y] circ U^G(x) circ (alpha_x)_* (u)}
\end{multline}
It suffices to check this for a basis element $u = [G_x/L] \in
U^G(\{\ast\}) = A(G_x)$  
for some subgroup $L \subseteq G_x$. Let $\pr\colon G/L \to G/G_x$ be the projection.
We get from the definitions
\begin{eqnarray}
\lefteqn{\ch^G(M)_{[y]} \circ U^G(x) \circ (\alpha_x)_* ([G/L])} \nonumber
\\ 
& \hspace{20mm} = & 
\ch^G(M)_{[y]}([x \circ \pr\colon G/L \to M]) \nonumber
\\
& \hspace{20mm}= & 
\sum_{\substack{W\!K_y\cdot \tau \in \\W\!K_y\backslash \mor(y,x\circ \pr)}}
|(W\!K_y)_{\tau}|^{-1}.
\label{ch^G(M)_[y] circ U^G(x) circ (alpha_x)_* ([G/L])}
\end{eqnarray}
We get a $W\!K_y$-map
$$q\colon \mor(y,x\circ \pr) \to \mor(y,x), \hspace{5mm} \tau \mapsto \pr
\circ\tau.$$ 
We can write
$$ \mor(y,x \circ \pr) ~ = ~ 
\coprod_{\substack{W\!K_y\cdot \sigma \in \\W\!K_y\backslash \mor(y,x)}}
W\!K_y \times_{(W\!K_y)_{\sigma}} q^{-1}(\sigma).$$
The $(W\!K_y)_{\sigma}$-set $q^{-1}(\sigma)$ is a finite disjoint
union of orbits 
$$q^{-1}(\sigma) ~ = ~ \coprod_{i \in I(\sigma)} (W\!K_y)_{\sigma}/A_i.$$
This implies
$$ \mor(y,x \circ \pr) ~ = ~ 
\coprod_{\substack{W\!K_y\cdot \sigma \in \\W\!K_y\backslash \mor(y,x)}} ~ 
\coprod_{i \in I(\sigma)} W\!K_y/A_i$$
and hence
\begin{eqnarray*}
\sum_{\substack{W\!K_y\cdot \tau \in \\W\!K_y\backslash \mor(y,x \circ \pr)}}
|(W\!K_y)_{\tau}|^{-1}
& = & 
\sum_{\substack{W\!K_y\cdot \sigma \in \\W\!K_y\backslash \mor(y,x)}}~ 
\sum_{i \in I(\sigma)} |A_i|^{-1}.
\end{eqnarray*}
Putting this into \eqref{ch^G(M)_[y] circ U^G(x) circ (alpha_x)_* ([G/L])} 
yields
\begin{eqnarray}
\lefteqn{\ch^G(M)_{[y]} \circ U^G(x) \circ (\alpha_x)_* ([G/L])} \nonumber
\\ 
& \hspace{20mm}= & 
\sum_{\substack{W\!K_y\cdot \sigma \in \\W\!K_y\backslash \mor(y,x)}}~ 
\sum_{i \in I(\sigma)} |A_i|^{-1}.
\label{ch^G(M)_[y] circ U^G(x) circ (alpha_x)_* ([G/L]) second}
\end{eqnarray}
Obviously
\begin{eqnarray}
|q^{-1}(\sigma)| = |(W\!K_y)_{\sigma}| \cdot \sum_{i \in I(\sigma)} |A_i|^{-1}.
\label{|p^{-1}(sigma)|}
\end{eqnarray}
If $\sigma\colon  G/K \to G/G_x$ sends $g'K$ to $g'gG_x$ for appropriate $g \in G$ with
$g^{-1}Kg \subseteq G_x$, then there is a bijection
$$G_x/L^{g^{-1}Kg} \xrightarrow{\cong} q^{-1}(\sigma)$$
which maps $g'L$ to the morphism $R_{gg'}\colon G/K \to G/L, \hspace{2mm} g''K \to g''gg'L$. This
implies
\begin{eqnarray}
|q^{-1}(\sigma)| & = & |G_x/L^{g^{-1}Kg}| ~ = ~
 \ch_0^{G_x}([G_x/L])_{(K_{\sigma})}. 
\label{|p^{-1}(sigma)| second}
\end{eqnarray}
We conclude from the equations 
\eqref{ch^G(M)_[y] circ U^G(x) circ (alpha_x)_* ([G/L]) second},
\eqref{|p^{-1}(sigma)|} and 
\eqref{|p^{-1}(sigma)| second} above
that \eqref{ch^G(M)_[y] circ U^G(x) circ (alpha_x)_* (u)} follows for
$u = [G_x/L]$ and hence for all $u \in U^{G_x}(\{\ast\}) = A(G_x)$.

Now we are ready to prove Lemma \ref{character map and local Lefschetz class}.
We get from the Definition \ref{def: Lambda^G_{loc}(f)}
of $\Lambda_{\loc}^G(f)$  and from
\eqref{ch^G(M)_[y] circ U^G(x) circ (alpha_x)_* (u)}
\begin{multline}
\ch^G(M)_{[y]}(\Lambda_{\loc}^G(v)) 
~ = ~
\sum_{Gx \in G\backslash\Fix(f)}~  
\sum_{\substack{W\!K_y\cdot \sigma \in \\W\!K_y\backslash \mor(y,x)}} 
\\
|(W\!K_y)_{\sigma}|^{-1} \cdot \ch^{G_x}_0\left(\Deg_0^{G_x}\left((\id - T_xf)^c)\right)\right)_{(K_{\sigma})}.
\label{ch^G(M)_{[y](Lambda_{loc}^G(v))}}
\end{multline}
We conclude from
\eqref{Lambda^G and the character map}
\begin{eqnarray}
\ch^{G_x}_0\left(\Deg_0^{G_x}\left((\id -
T_xf)^c)\right)\right)_{(K_{\sigma})} 
& = & 
\deg\left(\left(\id_{T_xM^{K_{\sigma}}} -T_xf^{K_{\sigma}}\right)^c\right).
\label{ch^{G_x}_0((Deg^{G_x}((id - T_xf)^c))))_{(K_{sigma})}}
\end{eqnarray}
If $\sigma\colon G/K \to G/G_x$ is of the form $g'K \mapsto g'gG_x$
for appropriate 
$g \in G$ with $g^{-1}Kg \subseteq G_x$, then
\begin{eqnarray}
\deg\left(\left(\id_{T_xM^{K_{\sigma}}} -T_xf^{K_{\sigma}}\right)^c\right) 
~ = ~
\deg\left(\left(\id_{T_{gx}M^K(y)} - T_{gx}(f|_{M^K(y)})\right)^c\right)
\label{deg((id_{T_xM^{K_{sigma}}} -T_xf^{K_{sigma}}))^c)}
\end{eqnarray}
since $f$ is $G$-equivariant. There is a bijection of $W\!K_y$-sets
$$\coprod_{Gx \in G\backslash\Fix(f)} ~  
\mor(y,x) ~ 
\xrightarrow{\cong}  \Fix(f|_{M^K(y)})$$
which sends $\sigma \in \mor(y,x)$ to $\sigma(1K)\cdot x$.
Now Lemma \ref{character map and local Lefschetz class} follows from
\eqref{ch^G(M)_{[y](Lambda_{loc}^G(v))}}, 
\eqref{ch^{G_x}_0((Deg^{G_x}((id - T_xf)^c))))_{(K_{sigma})}} and
\eqref{deg((id_{T_xM^{K_{sigma}}} -T_xf^{K_{sigma}}))^c)}. \qed

Because of Lemma  \ref{lem: injectivity of the character map chi^G},
Lemma \ref{character map and Lefschetz number}
and Lemma \ref{character map and local Lefschetz class} the
equivariant Lefschetz fixed point 
Theorem \ref{the: equivariant Lefschetz theorem} follows from
the orbifold Lefschetz fixed point
Theorem \ref{the: the orbifold Lefschetz fixed point theorem}.\qed

%%%%%%%%%%%%%%%%%%%%%%%%%%   Section 6  %%%%%%%%%%%%%%%%%%%%%%%%%%%%%%%%%%%%%%%%

\typeout{-----------------------  Section 6  ------------------------}

\tit{Euler characteristic and index of a vector field in the
equivariant setting} 
{Euler characteristic and index of a vector field in the equivariant setting}

\begin{definition} 
\label{def: equivariant universal Euler characteristic}
Let $X$ be a finite $G$-$CW$-complex $X$. We define the
\emph{universal equivariant Euler characteristic} of $X$
\begin{eqnarray*}
\chi^G(X) & \in & U^G(X)
\end{eqnarray*}
by assigning to $[x\colon G/H \to X] \in \Is \Pi_0(G,X)$ the
{\lp}ordinary{\rp} Euler characteristic of the pair of finite $CW$-complexes
$(W\!H_x\backslash X^H(x),W\!H_x\backslash X^{>H}(x))$.
If $X$ is proper, we define its \emph{orbifold Euler characteristic}
\begin{eqnarray*}
\chi^{\qq G}(X) & := & \sum_{p \ge 0} ~ (-1)^p
\sum_{G\cdot e \in G\backslash
I_p(X)} |G_e|^{-1} 
\hspace{5mm} \in \qq,
\end{eqnarray*}
where $I_p(X)$ is the set of the open cells of the CW-complex $X$ {\lp}after
forgetting the group action{\rp} and $G_e$ is the isotropy group of $e$
under the $G$-action on $I_p(X)$.  
\end{definition}

One easily checks that
$\chi^G(X)$ agrees with the equivariant Lefschetz class $\Lambda^G(\id_X)$
(see Definition \ref{def: equivariant Lefschetz class Lambda^G(f)}). It
can also be expressed by counting equivariant cells. If $x\colon G/H \to X$
is an object in $\Pi_0(G,X)$ and 
$\sharp_p([x])$ is the number of equivariant $p$-dimensional cells of orbit 
type $G/H$ which meet the component $X^H(x)$, then
\begin{eqnarray}
\chi^G(X)([x]) & = & \sum_{p \ge 0} (-1)^p \cdot \sharp_p([x]).
\label{chi^G(X) = counting cells}
\end{eqnarray}
The universal equivariant Euler characteristic is the universal additive
invariant for finite 
$G$-$CW$-complexes in the sense of \cite[Theorem 6.7]{Lueck(1989)}.

The orbifold Euler characteristic $\chi^{\qq G}(X)$ can be identified with
$L^{\qq G}(\id_X)$ (see Definition \ref{def: orbifold Lefschetz
number}) or with the more general notion of the $L^2$-Euler characteristic
$\chi^{(2)}(X;\caln(G))$. 
In analogy with $L^{(2)}(f,f_0;\caln(G))$ (see Remark \ref{rem: L^G(f)}),
one can compute $\chi^{(2)}(X;\caln(G))$ in terms of $L^2$-homology
$$\chi^{(2)}(X;\caln(G)) ~ = ~ 
\sum_{p \ge 0} (-1)^p \cdot
\dim_{\caln(G)}\left(H_p^{(2)}(X;\caln(G)\right),$$ 
where $\dim_{\caln(G)}$ denotes the von Neumann dimension (see for instance
\cite[Section 6.6]{Lueck(2002)}).

We conclude from Lemma \ref{character map and Lefschetz number}:

\begin{lemma} \label{character map and Euler characteristic}
Let $X$ be a finite proper $G$-$CW$-complex $X$. Let
$[y]$ be an isomorphism class of objects $y\colon G/K \to X$ in $\pi_0(G,X)$.
Then
$$\ch^G(X)_{[y]}(\chi^G(X)) ~ = ~ \chi^{\qq W\!K_y}(X^K(y)).$$
\end{lemma}

Next we want to express the universal equivariant Euler characteristic
of a cocompact proper $G$-manifold $M$, possibly with boundary,
in terms of the zeros of an equivariant vector field.

Consider an \emph{equivariant vector field}
$\Xi$ on $M$, i.e., a $G$-equivariant section of the tangent bundle $TM$
of $M$. Suppose that $\Xi$ is \emph{transverse to the
zero-section} $i\colon M \to TM$, i.e., if $\Xi(x) = 0$, then
$T_{\Xi(x)}TM$ is the sum of the subspaces given by the images of
$T_x\Xi\colon T_xM \to T_{\Xi(x)}(TM)$ and 
$T_xi\colon T_xM \to T_{\Xi(x)}(TM)$. Any equivariant vector field on $M$ 
can be changed by an arbitrary small perturbation into one 
which is transverse to the zero-section. We want to assign to such a
$\Xi$ its \emph{equivariant index} as follows. 

Since $\Xi$ is transverse to the zero-section, the set $\Zero(\Xi)$ of
points $x \in M$ with $\Xi(x) = 0$ is discrete. Hence $G\backslash
\Zero(\Xi)$ is 
finite, since $G$ acts properly on $M$ and $G\backslash M$ is compact by
assumption. Fix $x \in \Zero(\Xi)$. The zero-section $i\colon M \to TM$
and the inclusion 
$j_x\colon T_xM \to TM$ induce an isomorphism of $G_x$-representations
$$T_xi \oplus T_xj_x\colon T_xM \oplus T_xM \xrightarrow{\cong}  
        T_{i(x)}(TM)$$
if we identify $T_{i(x)}(T_xM) = T_xM$ in the obvious way.  If $\pr_k$
denotes the projection onto the $k$-th factor for $k = 1,2$ we obtain a
linear $G_x$-equivariant isomorphism
\begin{eqnarray}
d_x\Xi\colon & 
T_xM ~ \xrightarrow{T_x\Xi} T_{\Xi(x)}(TM) ~
~ \xrightarrow{(T_xi \oplus T_xj_x)^{-1}} ~
T_xM \oplus T_xM ~\xrightarrow{\pr_2} ~ T_xM.~ &
\label{definition of d_xv}
\end{eqnarray}
Notice that we obtain the identity if we replace $\pr_2$ by $\pr_1$ in
the expression \eqref{definition of d_xv} above. 
The $G_x$-map $d_x\Xi$ induces a $G_x$-map
$(d_x\Xi)^c\colon T_xM^c \to T_xM^c$
on the one-point compactification. 
Define analogously to the \emph{local equivariant Lefschetz class}
(see Definition  \ref{def: Lambda^G_{loc}(f)})
the \emph{equivariant index} of $\Xi$ 
\begin{eqnarray}
i^G(\Xi) & \in & U^G(M) \label{definition of equivariant index i^G(v)}
\end{eqnarray}
by
$$i^G(\Xi) ~ := ~  \sum_{Gx \in G\backslash\Zero(\Xi)}~ 
U^G(x) \circ (\alpha_x)\left(\Deg_0^{G_x}((d_x\Xi)^c)\right).$$
We say that $\Xi$ points \emph{outward  at the boundary} if for each 
$x \in \partial M$ the tangent vector $\Xi(x) \in T_xM$ does not lie in
the subspace 
$T_x\partial M$ and is contained in the half space $T_xM^+$ of tangent vectors
$u \in T_xM$ for which there is a path $w\colon [-1,0] \to M$ with
$w(0) = x$ and $w'(0) = u$. 
If $\partial M = \emptyset$, this condition is always satisfied.

\begin{theorem}[Equivariant Euler characteristic and vector fields]
\label{the: chi^G(X) = i^G(v)}
Let\linebreak 
$M$ be a cocompact proper smooth $G$-manifold. Let $\Xi$ be a $G$-equivariant
vector field which is transverse to the zero-section and points outward
at the boundary. Then we get in $U^G(M)$
$$\chi^G(M) ~ = ~ i^G(\Xi).$$
\end{theorem}
\proof
Let $\Phi\colon  M \times (-\infty,0]  \to M$ be the flow
associated to the vector field $\Xi$. It is defined on
$M \times (-\infty,0]$ since $\Xi$ is equivariant and points outward and
$M$ is cocompact. 
Moreover, each map $\Phi_{-\epsilon}\colon M \to M$ is a
$G$-diffeomorphism and 
$G$-homotopic to $\id_M$ for $\epsilon > 0$. This implies
$$\Lambda^G(\Phi_{-\epsilon}) ~ = ~ \Lambda^G(\id_M) ~ = ~ \chi^G(M).$$
Because of the equivariant Lefschetz fixed point Theorem 
\ref{the: equivariant Lefschetz theorem} it remains to prove for some
$\epsilon > 0$ 
$$i^G(\Xi) ~ = ~ \Lambda_{\loc}^G(\Phi_{-\epsilon}).$$
If we choose $\epsilon > 0$ small enough, 
the diffeomorphism $\phi_{-\epsilon}\colon  M \to M$ will have as set
of fixed points  $\Fix(\phi_{-\epsilon})$ precisely
$\Zero(\Xi)$. It suffices to prove for $x \in \Zero(\Xi) =
\Fix(\Phi_{-\epsilon})$ 
$$\Deg_0^{G_x}((d_x\Xi)^c) ~ = ~ \Deg_0^{G_x}((\id -
T_x\Phi_{-\epsilon})^c).$$ 
Recall that the character map $\ch^{G_x}_0\colon A(G_x) \to \prod_{(H)
\in \consub(G_x)} \zz$  
of \eqref{character map ch^K_0 of A(K) for finite K} 
is injective. We conclude from \eqref{Lambda^G and the character map}
that  it suffices to prove for any subgroup $H \subseteq G_x$ 
the equality of degrees of self-maps of the closed orientable manifold
$((T_xM)^c)^H = ((T_xM)^H)^c$
\begin{eqnarray}
\deg\left(((d_x\Xi )^c)^H\right) 
& = & 
\deg\left(((\id - T_x\Phi_{-\epsilon})^c)^H\right).
\label{comparing the degrees of v and Phi}
\end{eqnarray}
It suffices to treat the case $H = \{1\}$ and $\dim(M) \ge 1$ --- the
other cases are  
completely analogous or follow directly from the definitions.
Since \eqref{comparing the degrees of v and Phi}  is of local nature, 
we may assume $M = \rr^n$ and $x = 0$.
In the sequel we use the standard identification $T\rr^n = \rr^n
\times \rr^n$. 
Then the vector field $\Xi$ becomes a smooth map $\Xi \colon  \rr^n
\to \rr^n$ with $\Xi(0) = 0$ and 
$d_0\Xi$ becomes the differential $T_0\Xi$. Let $\Phi$ be the flow
associated to $\Xi$. Choose $\epsilon\ge 0$ and  
an open neighborhood $U \subseteq \rr^n$ of $0$ 
such that $\Phi$ is defined on $U \times [-\epsilon,0]$. By Taylor's
theorem we can 
find a smooth map $\eta\colon U \times [-\epsilon,0] \to \rr^n$ such that
for $t \in [-\epsilon,0]$ and $u \in U$,
$$\phi_t(u) = u + t \cdot \Xi(u) + t^2 \cdot \eta_t(u)~.$$
This implies
$$T_0\phi_t = \id  + t \cdot \left(T_0 \Xi + t \cdot T_0\eta_t\right).$$
Since $[-\epsilon,0]$ is compact, we can find 
a constant $C$ independent of $t$ such that the operator norm of
$T_0\eta_t$ satisfies 
$\Vert T_0\eta_t\Vert < C$ for $t \in [-\epsilon,0]$. The differential
$T_0 \Xi$ is an isomorphism by assumption. 
Hence $T_0\Xi + t \cdot T_0\eta_t$ is invertible for $t \in [-D,0]$ if we put
$D: = \min\{\epsilon, C^{-1} \cdot ||(T_0\Xi)^{-1}||^{-1}\}$.
Hence we get for $t \in [-D,0]$ that $\id - T_0\Phi_t$ is invertible and
$$
\frac{\det(\id - T_0\Phi_t)}{|\det(\id - T_0\Phi_t)|} 
~ = ~ 
\frac{\det(T_0 \Xi + t \cdot T_0\eta_t)}{|\det(T_0 \Xi + t \cdot T_0\eta_t)|}
~ = ~ 
\frac{\det(T_0 \Xi)}{|\det(T_0 \Xi)|}.$$
Hence \eqref{comparing the degrees of v and Phi} follows if we
take $\epsilon > 0$ small enough.
This finishes the proof that Theorem \ref{the: chi^G(X) = i^G(v)}
follows from Theorem \ref{the: equivariant Lefschetz theorem}. \qed

\begin{remark} \label{rem: vanishing of chi^G(M) and non-vanishing
vector fields} 
\em Let $M$ be a proper cocompact $G$-manifold without boundary. 
If $M$ possesses a nowhere-vanishing equivariant vector field, then
$\chi^G(M) = 0$ by the Theorem
\ref{the: chi^G(X) = i^G(v)}. The converse is true if 
$M$ satisfies the weak gap hypothesis that 
$\dim(M^{>G_x}(x)) \le \dim(M^{G_x}(x)) - 2$ holds for each $x \in G_x$. 
The proof is done by induction over the orbit bundles and the
induction step is reduced  
to the non-equivariant case. The weak gap condition ensures
that $M^{G_x}(x) - M^{>G_x}(x)$ is connected for $x \in M$. It is satisfied
if all isotropy groups of $M$ have odd order.  For finite groups $G$
more information about this question  can be found in \cite[Remark 2.5
(iii) on page 32]{Waner-Wu(1988)}. 
\em
\end{remark}

\begin{example} \label{exa: infinite dihedral group}
\em Let $D$ be the infinite dihedral group $D = \zz \rtimes \zz/2 =
\zz/2 \ast \zz/2$. 
We use the presentation $D = \langle s,t \mid s^2 = 1, s^{-1}ts =
t^{-1}\rangle$ in the 
sequel. The subgroups $H_0 = \langle t \rangle$ and $H_1 = \langle
ts\rangle$ have order 
two and $\{\{1\},H_0,H_1\}$ is a complete system of representatives for
the conjugacy classes 
of finite subgroups of $D$. The infinite dihedral group $D$ acts on $\rr$
by $s \cdot r = -r$ and $t \cdot r = r + 1$ for $r \in \rr$. The
interval $[0,1/2]$ is a 
fundamental domain for the $D$-action.  There is a $D$-$CW$-structure
on $\rr$ with $\{n \mid n \in \zz\} \coprod \{n + 1/2\mid n \in \zz\}$
as zero-skeleton. 
Let $x_i\colon D/H_i \to \rr$ be the $D$-map sending $1H_i$ to $0$ for
$i = 0$ and to  
$1/2$ for $i = 1$. Let $y \colon D \to \rr$ be the $D$-map sending $1$ to $0$.
Then $\rr$ has two equivariant $0$-cells, which have $x_0$ and $x_1$ as
their characteristic maps, 
and one equivariant $1$-cell of orbit type $D/\{1\}$.

We get $\Is\Pi_0(D;\rr) = \{[x_0],[x_1],[y]\}$. Recall that
$U^D(\rr)$ is the free $\zz$-module with basis
$\Is\Pi_0(D;\rr)$. Hence we write 
$$U^D(\rr) = \zz\langle [x_0] \rangle \oplus \zz\langle [x_1] \rangle
\oplus \zz\langle [y] \rangle.$$ 
We conclude from \eqref{chi^G(X) = counting cells}
$$\chi^D(\rr) = [x_0] + [x_1] - [y].$$
This is consistent with the original Definition 
\ref{def: equivariant universal Euler characteristic}, since
\begin{eqnarray*}
\chi(W_DH_i\backslash(\rr^{H_i}(x_i),\rr^{>H_i}(x_i))) & = &
\chi(\{\text{pt.}\}) = 1; 
\\
\chi(W_D\{1\}\backslash(\rr^{\{1\}}(y), \rr^{>\{1\}}(y)) & = &
\chi([0,1/2],\{0,1/2\}) = -1. 
\end{eqnarray*}
The character map \eqref{definition of character map chi^G(X)} is given by
\begin{multline*}
\chi^D(\rr) ~ = ~ \left\{\begin{array}{ccc} 1 & 0 & 0 \\ 0 & 1 & 0 \\
  1/2 & 1/2 & 1 
  \end{array}\right\} ~ \colon ~
U^D(\rr) = \zz\langle [x_0] \rangle \oplus \zz\langle [x_1] \rangle
  \oplus \zz\langle [y] \rangle 
\\
~ \to ~  \zz\langle [x_0] \rangle \oplus \zz\langle [x_1] \rangle
  \oplus \zz\langle [y] 
\rangle
\end{multline*}
since the $D$-set $\mor(y,x_i)$ is $D/H_i$ and the sets $\mor(x_i,y)$
and $\mor(x_i,x_j)$ 
for $i \not= j$ are empty. The character map sends $\chi^D(\rr)$ to
the various orbifold 
Euler characteristics which therefore must be
\begin{eqnarray*}
\chi^{\qq W_DH_i}(\rr^{H_i}(x_i)) & = & 1;
\\
\chi^{\qq D}(\rr) & = & 0.
\end{eqnarray*}
One easily checks that this is consistent with the Definition 
\ref{def: equivariant universal Euler characteristic} of the orbifold Euler
characteristics.

Let $\Xi$ be a vector field on $\rr$. Under the standard
identification $T\rr = \rr \times 
\rr$ this is the same as a function $\Xi \colon \rr \to \rr$. The
vector field $\Xi$ is 
transverse to the zero-section if and only if the function $\Xi$ satisfies
$\Xi(z) = 0 \Rightarrow \Xi'(z) \not= 0$. The vector field $\Xi$ is
$D$-equivariant if and only if 
$\Xi (-z) = -\Xi (z)$ and $\Xi (z) = \Xi (z+1)$ holds for all $z \in
\rr$. Let $\Xi$ be a 
$D$-equivariant vector field transverse to the zero-section. For
example, we can take $\Xi (z)=\sin (2\pi z)$.

We conclude $\Xi (0) = 0$ from $\Xi (-z) = -\Xi (z)$ and $\Xi(1/2) = 0$ from
$\Xi(1-z) = -\Xi(z)$. Let $z_0 = 0 < z_1 < z_2 < \ldots < z_r = 1/2$ 
be the points $z \in  [0,1/2]$ for which $\Xi (z) = 0$. (The example
of $\sin (2\pi z)$ shows the minimum value of $r$ is $1$.)
For $i \in \{0,1, \ldots  , r\}$ put $\delta_i =
\frac{\Xi'(z_i)}{|\Xi'(z_i)|}$. 
Since $\Xi$ is different from zero on $(z_i,z_{i+1})$, we have
$\delta_{i+1} = - \delta_i$ 
for $i \in \{0,1, \ldots  , r-1\}$. Let $d_{z_i}\Xi \colon T_{z_i}\rr
\to T_{z_i}\rr$ be the 
map associated to $\Xi$ at $z_i$ (see \eqref{definition of d_xv}); this
is simply multiplication by $\Xi'(z_i)$ under the
standard identification $T_{z_i}\rr = \rr$. The degree of the map
$(d_{z_i}\Xi)^c$ induced on 
$\rr^c = S^1$ is $\delta_i$. For $i = 0,r$ the isotropy group 
$ D_{ z_i}= \zz/2$ acts on $\rr$ by $- \id$. Hence the degree
$\deg\left((d_{z_i}\Xi )^{D_{ z_i}})^c\right)$ is by definition $1$
since $\dim(\rr^{D_{z_i}}) = 0$. 
We conclude from \eqref{Lambda^G and the character map} for $i = 0,r$
\begin{eqnarray*}
\Deg^{D_{z_i}}(d_{z_i}\Xi) & = & [H/H] + \frac{-1 + \delta_i}{2} \cdot
[H/\{1\}] 
\hspace{5mm} \in A(D_{z_i})
\end{eqnarray*}
and for $i \in \{1,2, \ldots , r-1\}$
\begin{eqnarray*}
\Deg^{D_{z_i}}(d_{z_i}\Xi) & = & \delta_i \cdot[\{1\}]
\hspace{5mm} \in A(\{1\}).
\end{eqnarray*}
Hence we get for the equivariant index
\begin{multline*}
i^D(\Xi) ~ = ~ 
[x_0] + \frac{-1 + \delta_0}{2} \cdot [y] + [x_1] + \frac{-1 +
\delta_r}{2} \cdot [y] + 
\sum_{i=1}^{r-1} \delta_i \cdot [y]
\\
= ~ [x_0] + [x_1] + \left(-1 + \frac{\delta_0}{2} + \frac{\delta_r}{2}
+ \sum_{i = 1}^{r-1} 
  \delta_i\right) \cdot [y]
\\
= ~ [x_0] + [x_1] - [y] + \delta_0 \cdot \left(\frac{1}{2} +
\frac{(-1)^r}{2}   + \sum_{i = 1}^{r-1} 
  (-1)^i \right) \cdot [y_0] 
\\
= ~ [x_0] + [x_1] - [y] + \delta_0 \cdot 0  \cdot [y_0] 
\\
= ~ [x_0] + [x_1] - [y].
\end{multline*}
This is consistent with Theorem \ref{the: chi^G(X) = i^G(v)}. 
\em
\end{example}

%%%%%%%%%%%%%%%%%%%%%%%%%%%%%%%%  Section 7 %%%%%%%%%%%%%%%%%%%%%%%%%%%%%

\typeout{-----------------------  Section 7  ------------------------}

\tit{Constructing equivariant manifolds with given component structure
and universal equivariant Euler characteristic}
{Constructing equivariant manifolds with given component structure and
universal equivariant Euler characteristic}

In this section we discuss the problem of whether 
there exists a proper smooth $G$-manifold $M$
with prescribed sets $\pi_0(M^H)$ for $H \subseteq G$, and whether 
$\chi^G(M)$  can realize a given element
in $U^G$. 

Let $\Orcat$ be the orbit category for the family of finite subgroups,
i.e., objects are homogeneous spaces $G/H$ for finite subgroups
$H\subseteq G$ and 
morphisms are $G$-maps. A contravariant $\Orcat$-set $S$ is a
contravariant functor 
from $\Orcat$ to the category of sets. Given a $G$-space $X$, define
a contravariant $\Orcat$-set $\underline{\pi_0(X)}$ by sending
$G/H$ to $\pi_0(X^H) = \pi_0(\map_G(G/H,X))$.
Next we investigate under which conditions a 
contravariant $\Orcat$-set $S$ arises
from a finite proper $G$-$CW$-complex.

\begin{lemma} \label{lem: realization of underline{pi_0(X)}}
For a contravariant $\Or(G;\calfin)$-set 
$S\colon \Or(G;\calfin) \to \Sets$ the following assertions are equivalent:

\begin{enumerate}

\item \label{lem: realization of underline{pi_0(X)}: set}
There are only finitely many elements $(H) \in \consub(H)$ with
$S(G/H)\not= \emptyset$. For any finite subgroup $H\subseteq G$ the set
$W\!H\backslash S(G/H)$ is finite and the isotropy group $W\!H_s$ of
each element 
$s \in S(G/H)$ is finitely generated;

\item \label{lem: realization of underline{pi_0(X)}: CW}
There is a proper finite $G$-$CW$-complex $X$ such that there exists
a natural equivalence $T\colon \underline{\pi_0(X)} \xrightarrow{\cong} S$.

\end{enumerate}
\end{lemma}
\proof
\eqref{lem: realization of underline{pi_0(X)}: CW} $\Rightarrow$
\eqref{lem: realization of underline{pi_0(X)}: set}
Since $X$ is finite, there are finitely many elements $(K_1)$, $(K_2)$,
$\ldots$, $(K_m)$ in $\{(K) \in \consub(G)\mid |K| < \infty\}$ 
such that for each equivariant cell $G/H \times D^n$
there is $i \in \{1,2, \ldots, m\}$ with $(H) = (K_i)$. Hence any
subgroup $H \subseteq G$ with $X^H \not= \emptyset$ is conjugate to
a subgroup of one of the $K_i$'s. Since a finite group has only
finitely many subgroups, 
the set $\{(H) \in \consub(G)\mid X^H \not= \emptyset\}$ is finite. 

Since $X$ is finite and proper and $W\!H\backslash(G/K^H)$ is finite for 
each finite group $H \subseteq G$ and subgroup $K \subseteq G$, 
the $W\!H$-$CW$-complex $X^H$ is finite proper for each finite
subgroup $H \subseteq G$. 
Hence  the quotient space $W\!H\backslash X^H$ is a finite
$CW$-complex and has only 
finitely many components. Since $W\!H\backslash\pi_0(X^H) \cong
\pi_0(W\!H\backslash 
X^H)$, the set $W\!H\backslash\pi_0(X^H)$ is finite.

Consider a finite subgroup $H\subseteq G$ and a component $C \in
\pi_0(X^H)$. Let 
$W\!H_C$ be the isotropy group $C$. Then $C$ 
is a connected proper finite $W\!H_C$-$CW$-complex. The long exact
homotopy sequence 
of the fibration 
\[
C \to EW\!H_C \times_{W\!H_C} C \to BW\!H_C
\]
shows that $W\!H_C$ is a quotient of $\pi_1(EW\!H_C \times_{W\!H_C}
C)$. Since for any finite subgroup $K \subseteq W\!H_C$ the 
$W\!H_C$-space $EW\!H_C \times_{W\!H_C} W\!H_C/K$ is homotopy equivalent to 
$BK$ and hence to a $CW$-complex of finite type (a $CW$-complex 
for which each skeleton is finite), $C$ is built  
out of finitely many cells $G/K \times D^i$, $K \subseteq G$ finite, 
and $EW\!H_C \times_{W\!H_C} W\!H_C/K$ has the homotopy type of
$CW$-complex of finite type.  
This implies that $\pi_1(EW\!H_C \times_{W\!H_C} C)$ is finitely
generated. Hence 
$W\!H_C$ is finitely generated.
\\[2mm]
\eqref{lem: realization of underline{pi_0(X)}: set} $\Rightarrow$
\eqref{lem: realization of underline{pi_0(X)}: CW}
Choose an ordering $(H_1), (H_2), \ldots, (H_r)$ on the set
$\{(H) \in \consub(G) \mid S(G/H) \not= \emptyset\}$ such that
$H_i$ is subconjugate to $H_j$ only if $i \ge j$ holds.
Define 
$$X_0 ~ := ~ \coprod_{i=1}^s \coprod_{W\!H_i\backslash S(G/H_i)} G/H_i.$$
Define a transformation $\phi_0\colon \underline{\pi_0(X_0)} \to S$ as
follows. 
Fix an object $G/K \in \Or(G;\calfin)$. 
Then  $\underline{\pi_0(X_0)}$ evaluated at this object $G/K$
is given by
$$ \coprod_{i=1}^s \coprod_{W\!H_i\backslash S(G/H)} \map_G(G/K,G/H_i)$$
since $\pi_0(G/H_i^K) =  \map_G(G/K,G/H_i)$. Now require that
$\phi_0(G/K)$ sends 
a $G$-map $\sigma \in  \map_G(G/K,G/H_i)$ belonging to the summand
for $W\!H_i \cdot s \in W\!H_i\backslash S(G/H_i)$ to
$S(\sigma)(s)$.  One easily checks that
$\phi_0(G/H)\colon  \pi_0(X^H) \to S(G/H)$ is surjective for all $H
\subseteq G$. 

In the next step we attach equivariant one-cells to $X_0$ to get a
$G$-$CW$-complex $X$ together with a transformation
$\phi\colon \underline{\pi_0(X)} \to S$, such that the composition of
$\phi$   with the transformation $\underline{\pi_0(X_0)} \to
\underline{\pi_0(X)}$ 
induced by the inclusion is $\phi_0$, and $\phi(G/H)$ is bijective for
all $H \subseteq G$. We do this by constructing 
by induction a sequence of proper 
cocompact $G$-$CW$-complexes $X_0 \subseteq X_1 \subseteq \ldots
\subseteq X_r$, together with  transformations
$\phi_i\colon \underline{\pi_0(X_i)} \to S$ such that the composition of
$\phi_i$   with the transformation $\underline{\pi_0(X_{i-1})} \to
\underline{\pi_0(X_i)}$ 
induced by the inclusion is $\phi_{i-1}$, $\phi_i(G/H_j)$ is
bijective for all $j \le i$, and $X_i$ is obtained from $X_{i-1}$ by
attaching finitely many equivariant cells of the type $G/H_i 
\times D^1$. Then we can take $X = X_r$ and $\phi = \phi_r$.

The induction beginning is $X_0$ together with $\phi_0$. The 
induction step from $i-1$ to $i$ is done as follows.
Consider $W\!H_i \cdot s \in W\!H_i\backslash S(G/H_i)$. Let
$(W\!H_i)_s$ be the isotropy group of 
$s$ under the $W\!H_i$-action on $S(G/H_i)$. The preimage of $s$ under 
$\phi_{i-1}(G/H_i) \colon \pi_0(X_{i-1}^{H_i}) \to S(G/H_i)$ consists
of finitely many 
$(W\!H_i)_s$-orbits since $W\!H_i\backslash \pi_0(X_{i-1}^{H_i})$ is finite by
the implication
\eqref{lem: realization of underline{pi_0(X)}: CW} $\Rightarrow$
\eqref{lem: realization of underline{pi_0(X)}: set} which we have
already proved. 
Let $u(s)_1$, $u(s)_2$, $\ldots, u(s)_v$ be a system of generators of
$(W\!H_i)_s$ which contains the unit element $1 \in (W\!H_i)_s$. 
Fix  elements  $C(s)_1$, $C(s)_2$, $\ldots$, $C(s)_{n(s)}$ in 
$\phi_{i-1}^{-1}(G/H_i)(s)$ such that 
$$(W\!H_i)_s \backslash \phi_{i-1}(G/H_i)^{-1}(s) ~ = ~ 
\{(W\!H_i)_s \cdot C(s)_i \mid i = 1,2 \ldots , n(s)\}$$
and $(W\!H_i)_s \cdot C(s)_j = (W\!H_i)_s \cdot C(s)_k$ implies $j = k$. 
Now attach for each $W\!H_i \cdot s \in W\!H_i\backslash S(G/H_i)$, each 
generator $u(s)_i$, each $C(s)_j$ an equivariant cell $G/H_i \times D^1$ to $X_{i-1}$ 
such that $\{1H_i\} \times D^1$ connects $C_1(s)$ and $u(s)_i \cdot C_j(s)$.
The resulting $G$-$CW$-complex is the desired $G$-$CW$-complex $X_i$, one easily constructs
the desired transformation $\phi_i$ out of $\phi_{i-1}$.
\qed

Notice that for a finite group $G$ the statement
\eqref{lem: realization of underline{pi_0(X)}: set} in Lemma
\ref{lem: realization of underline{pi_0(X)}} is equivalent to the statement
that $S(G/H)$ is finite for all subgroups $H \subseteq G$. If
we take as $\Or(G;\calfin)$-set the functor $S$ which sends $G/\{1\}$ to $\{\ast\}$ and
$G/H$ for $H \not= \{1\}$ to $\emptyset$, then Lemma
\ref{lem: realization of underline{pi_0(X)}} boils down to the statement that 
there exists a connected finite free $G$-$CW$-complex $X$ if and only if 
$G$ is finitely generated.

Now given a contravariant $\Orcat$-set $S$, we define $U^G(S)$ to be the
free abelian group on $\coprod_{(H)\in\consub(G)}WH\!\backslash S(G/H)$.
If $S$ satisfies the equivalent conditions of Lemma
\ref{lem: realization of underline{pi_0(X)}}, then clearly this is
naturally isomorphic to $U^G(X)$, with $X$ as in part (b) of
Lemma \ref{lem: realization of underline{pi_0(X)}}.

Next we  prove that is $X$ is a finite proper $G$-$CW$-complex,
then any element $u\in U^G(X)$ can be realized  from $\chi^G$ of
a manifold. More precisely, there is a $G$-map $f\colon M\to X$ with
$M$ a $G$-manifold, such that $U^G(f)$ is an isomorphism sending
$\chi^G(M)$ to $u$.
We are grateful to Tammo tom Dieck for pointing out to us 
the use of the multiplicative induction
in the proof of the next result.

\begin{lemma} \label{lem: realization of elements in U^G(X) by chi^G(M)}
Let $X$ be a finite proper $G$-$CW$-complex and $u \in U^G(X)$.
Then there is a proper cocompact $G$-manifold $M$ without boundary together with
a $G$-map $f\colon M \to X$ with the following properties:

For any $x \in M$ the $G_x$-representation $T_xM$ is a multiple
of the regular $G_x$-representation $\rr[G_x]$ for $G_x$ the isotropy group of
$x \in X$. The dimensions of the components $C \in \pi_0(M)$ are all
equal. The components 
of $M^H$ are orientable manifolds for each $H \subseteq G$. The induced map
$\pi_0(f^H)\colon \pi_0(M^H) \to \pi_0(X^H)$ is bijective for each
finite subgroup $H \subseteq G$. The induced map
$$U^G(f)\colon U^G(M) \xrightarrow{\cong} U^G(X)$$ 
is bijective and sends $\chi^G(M)$ to $u$.
\end{lemma}
\proof In the first step we want to reduce the claim to the case,
where $X$ is  
a finite proper $1$-dimensional $G$-$CW$-complex such that for any
$x\colon G/H \to X$ 
there is a zero cell $G/H$ which  meets $X^H$.
 
If the $G$-map $f\colon X \to Y$ of proper $G$-$CW$-complexes
induces bijections $\pi_0(f^H)\colon \pi_0(X^H) \xrightarrow{\cong}
\pi_0(Y^H)$ for all 
finite subgroups $H \subseteq G$, then the induced map
$U^G(f)\colon U^G(X) \to U^G(Y)$ is bijective by 
\eqref{identifying Is Pi_0(G,X)}. In particular the inclusion of the $1$-skeleton
$i_1\colon X_1 \to X$ induces a bijection $U^G(i_1)$.

Since $X_1$ is finite and proper, $W\!H\backslash \pi_0(X^H)$ is finite
for all finite subgroups $H\subseteq G$ and the 
set $\{(H) \in \consub(G) \mid X^H \not= \emptyset\}$ is finite
(see Lemma \ref{lem: realization of underline{pi_0(X)}}).  We conclude
from \eqref{identifying Is Pi_0(G,X)} that $\Is\pi_0(G,X_1)$ is finite. Fix for any
$[x\colon G/H \to X_1]$ a representative $x\colon G/H \to X$ whose
image lies in $X_0$. 
This is possible by the equivariant cellular approximation theorem. Define a
finite proper $G$-$CW$-complex $Y$ by the $G$-pushout diagram
\comsquare{\coprod_{\substack{[x\colon G/H \to X] \in
\\\Is\pi_0(G,X)}}~  G/H \times \{0\}} 
{\coprod_{\substack{[x\colon G/H \to X] \in \\\Is\pi_0(G,X)}} ~ x }{X_1}
{i_2}{i_3}
{\coprod_{\substack{[x\colon G/H \to X] \in \\\Is\pi_0(G,X)}} ~ G/H \times [0,1]}{}{Y} 
where $i_2$ is the inclusion. Since $i_2$ is a $G$-homotopy equivalence,
$i_3$ is a $G$-homotopy equivalence. Let $i_3^{-1}\colon Y \to X_1$ be a $G$-homotopy inverse.
Then $U^G(i_1 \circ i_3^{-1})\colon U^G(Y) \to U^G(X)$ is a bijection and
$Y$ is a finite proper $G$-$CW$-complex such that 
for each $G$-map $y\colon G/H \to Y$ there is a zero cell $G/H$ which meets
$Y^H(y)$, namely $G/H \times \{1\} \subseteq G/H \times [0,1]$ for the corresponding
$[x] \in \Is \Pi_0(G,X)$ in the pushout diagram above.
Hence we can assume without loss of generality $X = Y$ in the sequel.

Fix a number $n$ such for any $H \subseteq G$ with $X^H \not= \emptyset$ the order
$|H|$ divides $n$.  Let $G/H_1$, $G/H_2$, $\ldots$, $G/H_r$ be the
equivariant zero-cells of $X$. Let $N_i$ be a $4n/|H_i|$-dimensional closed oriented
manifold. Let $\prod_{H_i} N_i$ be the closed $H_i$-manifold with the $H_i$-action coming from permuting
the factors. This is called the\emph{ multiplicative induction} or \emph{coinduction} of 
$N_i$. One easily checks that for $K \subseteq H_i$ 
the $K$-fixed point sets of $\prod_{H_i} N_i$ is diffeomorphic to $\prod_{k= 1}^{|H_i/K|} N_i$ and hence
a closed connected orientable manifold which is non-empty. Moreover,
the $(H_i)_x$-representation $T_x\left(\prod_{H_i} N_i\right)$ is a multiple of the regular real
$(H_i)_x$-representation for all $x \in \prod_{H_i} N_i$. The manifolds $N_i$ will be specified later.
Given an orthogonal $H$-representation $V$ of a finite group
$H$, we denote  by $DV$ and $SV$ the unit
disk and the unit sphere and  by $\inte(DV) = DV - SV$ the interior of $DV$.
Define the $4n$-dimensional proper cocompact $G$-manifold $M_0$ to be
$$M_0 ~ = ~ \coprod_{i=1}^r G \times_{H_i} \left(\prod_{H_i}  N_i\right).$$
Let $f_0\colon M_0 \to X$ be the map which is on 
$G \times_{H_i} \left(\prod_{H_i}  N_i\right)$ given by the canonical
projection  onto the cell $G/H_i = G \times_{H_i} \{\ast\}$. 

Fix a $G$-pushout describing, how $X = X_1$ is obtained from $X_0$ by
attaching equivariant one-cells
\comsquare
{\coprod_{j=1}^s G/K_j \times S^0}{q}{X_0}
{}{}
{\coprod_{j=1}^s G/K_j \times D^1}{Q}{X}

Consider a $G$-map $\sigma \colon  G/K \to G/H$ for finite subgroups
$H,K \subseteq G$. 
Suppose that $|H|$ and $|K|$ divide $n$. Choose $g \in G$ such that
$g^{-1}Kg \subseteq H$ and $\sigma$ sends $g'K$ to $g'gH$.
Let $c_g \colon K \to H$ be the injective group homomorphism
$g' \mapsto g^{-1}g'g$. The $K$-representations $\rr[K]^{4n/|K|}$ and 
$c_g^*\rr[H]^{4n/|H|}$, which is obtained from the $H$-representation
$\rr[H]^{4n/|H|}$ by restriction with $c_g$,  are isomorphic. 
Choose an isometric $c_g$-equivariant linear isomorphism
$$\phi\colon \rr[K]^{4n/|K|} \to \rr[H]^{4n/|H|}.$$
Choose a small number $\epsilon > 0$ and an element $w \in \rr[H]^{4n/|H|}$
such that $||w|| = 1$, the $H$-isotropy group of $w$ is
$g^{-1}Kg$ and the distance of two distinct points in the $H$-orbit
through $w$ is larger than $3\epsilon$.  The following map is a $G$-embedding
$$\psi\colon G \times_K D\rr[K]^{4n/|K|} \to G \times_H \rr[H]^{4n/|H|},
\hspace{5mm} (g',v) ~ \mapsto ~ (g'g,\epsilon \cdot \phi(v) + w)$$ 
such that the following diagram with the canonical projections as
vertical arrows commutes:
\comsquare{G \times_K D\rr[K]^{4n/|K|}}{\psi}{G \times_H \rr[H]^{4n/|H|}}
{\pr}{\pr}
{G/K}{\sigma}{G/H}
Using the construction above with appropriate choices of $w$ and $\epsilon$
we can find a $G$-embedding
$$\Psi\colon \coprod_{j=1}^s G\times_{K_j} \left(D\rr [K_j]^{4n/|K_j|}
\times S^0\right)  
\to  M_0$$
such that the following diagram commutes:
\comsquare{\coprod_{j=1}^s G\times_{K_j} \left(D\rr[K_j]^{4n/|K|}
\times S^0\right)} 
{\Psi}{M_0}
{\pr}{f_0}
{\coprod_{j=1}^s G/K_j \times S^0}{q}{X_0}
where $\pr$ is the obvious projection.  Let $M_0'$ be the proper
cocompact $G$-manifold 
with boundary
$$M_0' ~ = ~ 
M_0 - \Psi\left(\coprod_{j=1}^s G\times_{K_j} 
\left(\inte\left(D\rr [K_j]^{4n/|K_j|}\right) \times S^0\right)\right).$$
Define a smooth $G$-manifold $M$ by the following pushout
\comsquare{\coprod_{j=1}^s G\times_{K_j} \left(S\rr [K_j]^{4n/|K_j|}
\times S^0\right)} 
{\Psi|}{M_0'}
{}{}
{\coprod_{j=1}^s G\times_{K_j} \left(S\rr [K_j]^{4n/|K_j|} \times D^1\right)}
{\overline{\Psi|}}{M}
where $\Psi|$ denotes the restriction of $\Psi$. The map $f_0\colon
M_0 \to X_0$ defines a 
map $f_0'\colon M_0' \to X_0$ by restriction. 
Since $(S\rr [K_j]^{4n/|K|})^L$ for each $L \subseteq K_j$ and
$\left(\prod_{H_i} N_i\right)^L$ for each $L \subseteq H_i$ are
non-empty and connected, the maps 
$f_0'$ and the projections 
$\pr\colon \coprod_{j=1}^s G\times_{K_j} \left(S\rr [K_j]^{4n/|K|}
\times S^0\right) \to 
\coprod_{j=1}^s G/K_j \times S^0$ and
$\pr\colon \coprod_{j=1}^s G\times_{K_j} \left(S\rr [K_j]^{4n/|K|}
\times D^1\right) \to 
\coprod_{j=1}^s G/K_j \times D^1$  induce on the $L$-fixed point set
$1$-connected maps 
for each $L \subseteq G$. Hence the $G$-map $f\colon M \to X$, which
is induced by  
these three maps and the pushout property,
induces a $1$-connected map on the $L$-fixed point sets for each $L
\subseteq G$. 
In particular $\pi_0(f^L)\colon \pi_0(M^L) \to \pi_0(X^L)$ is
bijective for each 
$L \subseteq G$ and the map $U^G(f)\colon U^G(M) \to U^G(X)$ is bijective
by \eqref{identifying U^G(X)}.
Obviously the $G_x$-representation $T_xM$ for any point $x \in M$ is
a multiple of the regular $G_x$-representation $\rr[G_x]$.

The following diagram  commutes 
$$\begin{CD}
U^G(M_0) @< U^G(i_0) < \cong <  U^G(M_0') @> U^G(i_4) > > U^G(M)
\\
@V U^G(f_0) V\cong V @V U^G(f_0') V\cong V @V U^G(f) V\cong V
\\
U^G(X_0) @> \id > \cong >  U^G(X_0) @> U^G(i_5) >> U^G(X)
\end{CD}$$
where $i_0$, $i_4$ and $i_5$ denote the inclusions.

Next we compute $U^G(f)(\chi^G(M))$.
We get by the sum formula for the universal equivariant Euler characteristic
\cite[Theorem 5.4 on page 100]{Lueck(1989)}
\begin{eqnarray*}
\lefteqn{U^G(f)(\chi^G(M))} & & 
\\
&  =  &  U^G(i_5 \circ f_0)\left(\chi^G(M_0)\right)
\\
& & \hspace{6mm} - ~ 
U^G(i_5 \circ f_0 \circ \Psi) \left(\chi^G\left(\coprod_{j=1}^s G\times_{K_j} 
\left(\left(D\rr [K_j]^{4n/|K_j|}\right) \times S^0\right)\right) \right) 
\\
& &  \hspace{6mm} + ~ 
U^G(i_5 \circ f_0 \circ \Psi \circ i_6) \left(\chi^G\left(\coprod_{j=1}^s G\times_{K_j} 
\left(\left(S\rr [K_j]^{4n/|K_j|}\right) \times S^0\right)\right) \right) 
\\
& &  \hspace{6mm} - ~ 
U^G(i_5 \circ f_0' \circ \Psi|)\left(\chi^G\left(
\coprod_{j=1}^s G\times_{K_j} \left(S\rr [K_j]^{4n/|K_j|} \times S^0\right)\right)\right) 
\\
&  & \hspace{6mm} + ~
U^G(f \circ \overline{\Psi|})\left(\chi^G\left(
\coprod_{j=1}^s G\times_{K_j} \left(S\rr [K_j]^{4n/|K_j|} \times D^1\right)\right)\right), 
\end{eqnarray*} 
where $i_6$ is the inclusion.  Since 
$(S\rr [K_j]^{4n/|K_j|})^L$ is an odd-dimensional sphere and hence has
vanishing  (non-equivariant) Euler characteristic for each $L \subseteq
K_j$, the element $\chi^{K_j}\left(S\rr [K_j]^{4n/|K_j|}\right)$
in $U^{K_j}(\{\ast\}) = A(K_j)$ is sent to zero under the injective
map $\chi_0^{K_j}\colon A(K_j) \to \prod_{(L) \in \consub(K_j)} \zz$. This
implies 
$$\chi^{K_j}\left(S\rr [K_j]^{4n/|K_j|}\right) = 0$$
and hence 
\begin{eqnarray*}
\chi^G\left(\coprod_{j=1}^s G\times_{K_j} 
\left(S\rr [K_j]^{4n/|K_j|} \times S^0\right)\right) & = & 0.
\end{eqnarray*}
The space $D\rr [K_j]^{4n/|K_j|}$ is $K_j$-homotopy equivalent to
$\{\ast\}$. Hence
\begin{eqnarray*}
\lefteqn{U^G(i_5 \circ f_0 \circ \Psi)
\left(\chi^G\left(\coprod_{j=1}^s G\times_{K_j}  
\left(D\rr [K_j]^{4n/|K_j|} \times S^0\right)\right)\right)}
\\
& \hspace{5mm} = &  
U^G(i_5 \circ f_0 \circ \Psi \circ i_7)
\left(\chi^G\left(\coprod_{j=1}^s G\times_{K_j} S^0\right)\right)  
\\
& \hspace{5mm} = & 
2 \cdot \sum_{j=1}^s U^G(x_j)(\chi^G(G/K_j)),
\end{eqnarray*}
where $i_7$ is the inclusion and  $x_j^1\colon G/K_j \to X$ is the
restriction of the characteristic 
map of the one cell $G/K_j \times D^1$ to $G/K_j \times \{1/2\}$. If
$x_i^0\colon G/H_i \to X$ is the characteristic map of the
$0$-cell $G/H_i$ and $\alpha_i\colon H_i \to G$ the inclusion, we conclude
\begin{eqnarray}
U^G(f)(\chi^G(M)) & = & 
\sum_{i=1}^r U^G(x_i^0) \circ
(\alpha_i)_*\left(\chi^{H_i}\left(\prod_{H_i} N_i\right)\right)
\nonumber 
\\ & & \hspace{15mm} 
~  - ~ 2 \cdot \sum_{j=1}^s U^G(x_j)(\chi^G(G/K_j)).
\label{U^G(f)chi^G(M)): 1}
\end{eqnarray}
The element $\chi^{H_i}\left(\prod_{H_i} N_i\right) \in A(H_i)$ is
sent under the injective character map 
$$\ch^{H_i}_0\colon A(H_i) \to \prod_{(K) \consub(H_i)} \zz,
\hspace{5mm} [S] \mapsto |S^K|$$ 
to $(\chi(N_i)^{|H/K|} \mid (K) \in \consub(H_i)\}$. This implies
\begin{multline}
\chi^{H_i}\left(\prod_{H_i} N_i\right)   ~ = ~ \chi(N_i) \cdot [H_i/H_i] \\
+ \sum_{\substack{(K) \in \consub(H_i),\\ K \not= H_i}}
\lambda_{(K)}(\chi(N_i)) \cdot [H_i/K] 
\label{a(H_i) = [H_i/H_i] + higher terms}
\end{multline}
for appropriate functions $\lambda_{(K)}\colon \zz \to \zz$. 
We can order the equivariant zero-cells $G/H_i$
of $X_0$ such that $H_i$ is subconjugate to $H_j$ only if $i \ge j$
holds. We conclude from 
\eqref{a(H_i) = [H_i/H_i] + higher terms} that for an appropriate map
$\mu$, the composition  
$$
\bigoplus_{i=1}^r \zz \xrightarrow{\mu} \bigoplus_{i=1}^r \zz = U^G(X_0) \xrightarrow{U^G(i_8)} U^G(X_1)
$$
sends $\{\chi(N_i) \mid i = 1,2 \ldots , r\}$ to 
$\sum_{i=1}^r U^G(x_i^0) \circ (\alpha_i)_*\left(\chi^{H_i}\left(\prod_{H_i} N_i\right)\right)$,
where $i_8 \colon X_0 \to X$ is the inclusion and $\mu$ is given by
$$\mu(a_1, a_2, \ldots, a_r) ~ = ~ \left(\sum_{j=1}^r \mu_{1,j}(a_j), 
\sum_{j=1}^r \mu_{2,j}(a_j), \ldots , \sum_{j=1}^r \mu_{r,j}(a_j),\right)$$
for (not necessarily linear) maps $\mu_{i,j}\colon \zz \to \zz$ for which $\mu_{i,j} = 0$ for
$i > j$ and $\mu_{i,i} = \id$.  The map $U^G(i_8)$ is surjective and the map
$\mu $ bijective. Since for any integer $k$ and any positive integer $l$ there is a closed connected 
$4l$-dimensional manifold $N$ with $\chi(N) = k$, we can find
appropriate $N_i$ with the right Euler characteristics $\chi(N_i)$
such that for a given element  
$u \in U^G(X)$
\begin{multline}
\sum_{i=1}^r U^G(x_i^0) \circ
(\alpha_i)_*\left(\chi^{H_i}\left(\prod_{H_i} N_i\right)\right) 
~ = ~ u +  2 \cdot \sum_{j=1}^s U^G(x_j)(\chi^G(G/K_j)).
\label{U^G(f)chi^G(M)): 2}
\end{multline}
We conclude from 
\eqref{U^G(f)chi^G(M)): 1} and \eqref{U^G(f)chi^G(M)): 2} that
$$U^G(f)(\chi^G(M)) ~ = u.$$
This finishes the proof of 
Lemma \ref{lem: realization of elements in U^G(X) by chi^G(M)}. \qed

Putting Lemmas \ref{lem: realization of underline{pi_0(X)}} 
and \ref{lem: realization of elements in U^G(X) 
by chi^G(M)} together gives the following result on the realization
problem:

\begin{theorem}
\label{the: U^G realization theorem}
Let $S$ be a contravariant $\Orcat$-set, and suppose that
there are only finitely many elements $(H) \in \consub(H)$ with
$S(G/H)\not= \emptyset$. Also assume that
for any finite subgroup $H\subseteq G$, the set
$W\!H\backslash S(G/H)$ is finite, and that the isotropy group $W\!H_s$ of
each element $s \in S(G/H)$ is finitely generated. 
Let $u\in U^G(S)$. Then there is a proper cocompact $G$-manifold $M$ 
{\lp}without boundary{\rp} and there is a natural equivalence $T\colon
\underline{\pi_0(M)} \xrightarrow{\cong} S$ sending $\chi^G(M)$
to $u$.
\end{theorem}

%%%%%%%%%%%%%%%%%%%%%%%%%%%%%%  References %%%%%%%%%%%%%%%%%%%%%%%%%%%%%% 

\typeout{-------------------------- references --------------------------}

\noindent Addresses: \hfill \\[3mm]
Wolfgang L\"uck \hfill Jonathan Rosenberg\\
Institut f\"ur Mathematik und Informatik\hfill Department of Mathematics\\
Westf\"alische Wilhelms-Universtit\"at  \hfill University of Maryland\\
Einsteinstr. 62 \hfill College Park, MD 20742-4015\\
48149 M\"unster \hfill U.S.A\\
Germany \hfill \\
\texttt{
lueck@math.uni-muenster.de \hfill jmr@math.umd.edu\\
wwwmath.uni-muenster.de/math/u/lueck \hfill
www.math.umd.edu/\raisebox{-.6ex}{\symbol{"7E}}jmr}
\\[5mm]
%\date{Last edited: 16.7.02 or later\\last compiled: \today}

 \end{document}